\documentclass[11pt,a4paper]{article}

\usepackage{amsmath,amssymb}
\newcommand{\e}{\varepsilon}
\newcommand{\va}{\varphi}
\newcommand{\D}{\Delta}

\newcommand{\n}{\nabla}
\newcommand{\N}{\frac{N}{2}}
\newcommand{\NN}{\frac{N}{p}}

\newcommand{\p}{\partial}

\newcommand{\R}{\mathbb{R}}

\newtheorem{definition}{Definition}
\newtheorem{theorem}{Theorem}
\newtheorem{proposition}{Proposition}
\newtheorem{corollaire}{Corollary}

\newtheorem{remarka}{Remark}

\newtheorem{lemme}{Lemma}

\title{Local well-posedness results for density-dependent incompressible fluids}
\author{Boris Haspot \thanks{Karls Ruprecht Universit\"at Heidelberg, Institut for Applied Mathematics,
Im Neuenheimer Feld 294,
D-69120 Heildelberg, Germany.
Tel. 49(0)6221-54-6112 }}
\date{}
\begin{document}
\maketitle
\begin{abstract}
This paper is dedicated to the study of the initial value problem for density dependent
incompressible viscous fluids in $\R^{N}$ with $N\geq2$. We address
the question of well-posedness for {\it large} data having critical Besov regularity  and we aim at stating
well-posedness in functional spaces as close as possible to the ones imposed in the
incompressible Navier Stokes system by Cannone,
Meyer and Planchon in \cite{CMP} where $u_{0}\in B^{\NN-1}_{p,\infty}$ with $1\leq p<+\infty$.
This improves the analysis of \cite{Da3}, \cite{Da4} and \cite{AP} where $u_{0}$ is considered belonging to $B^{\NN-1}_{p,1}$ with $1\leq p<2N$.
Our result relies on a new a priori estimate for transport equation introduce by Bahouri, Chemin and Danchin in
\cite{BCD} when the velocity $u$
is not considered Lipschitz.
\end{abstract}
\section{Introduction}
In this paper, we are concerned with the following model of incompressible viscous fluid with variable density:
\begin{equation}
\begin{cases}
\begin{aligned}
&\p_{t}\rho+{\rm div}(\rho u)=0,\\
&\p_{t}(\rho u)+{\rm div}(\rho u\otimes u)-{\rm div}(2\mu(\rho)Du)+\n\Pi=\rho f,\\
&{\rm div}u=0,\\
&(\rho,u)_{/t=0}=(\rho_{0},u_{0}).
\end{aligned}
\end{cases}
\label{NSC}
\end{equation}
Here $u=u(t,x)\in\R^{N}$ stands for the velocity field and $\rho=\rho(t,x)\in\R^{+}$ is the density,
$Du=\frac{1}{2}(\n u
+^{t}\n u)$ is the strain tensor. We denote by  $\mu$ the viscosity coefficients of the fluid,
which is assumed to satisfy $\mu>0$.
The term $\n\Pi$ (namely the gradient of the pressure) may be seen as the Lagrange multiplier associated to
the constraint
${\rm div}u=0$.
We supplement the problem with initial condition $(\rho_{0},u_{0})$ and an outer force $f$.
Throughout the paper, we assume that the space variable $x\in\R^{N}$ or to the periodic
box ${\mathbb{T}}^{N}_{a}$ with period $a_{i}$, in the i-th direction. We restrict ourselves the case $N\geq2$.\\
The existence of global weak solution for (\ref{NSC}) under the assumption that $\rho_{0}\in L^{\infty}$
is nonnegative, that
${\rm div}u_{0}$, and that $\sqrt{\rho_{0}}u_{0}\in L^{2}$ has been studied by different authors. It
is based on the energy
equality:
\begin{equation}
\|\sqrt{\rho}u(t)\|_{L^{2}}^{2}+\int_{0}^{t}\|\sqrt{\mu(\rho)}Du(\tau)\|_{L^{2}}^{2}d\tau=\|\sqrt{\rho_{0}}
u_{0}\|_{L^{2}}^{2}
+\int 2(\rho f\cdot u)(\tau,x)d\tau\,dx.
\label{02}
\end{equation}
Using (\ref{02}) and the fact that the density is advected by the flow of $u$ so that the $L^{p}$ norms of $\rho$
are (at least formally) conserved during the evolution, it is then possible to use compactness methods to prove the existence
of global weak solution.
This approach has been introduced by J. Leray in 1934 in the homogeneous case (i.e $\rho=1$) and no external force.
For the non-homogeneous equation (\ref{NSC}), we refer to \cite{AKM} and to \cite{Li} for an overview of results
on weak solution. Some recent improvements have been obtained by B. Desjardins in \cite{De1},\cite{De2} and \cite{De3}.\\
The question of unique resolvability for (\ref{NSC}) has been first addressed by O. Ladyzhenskaya and V. Solonnikov
in the late
seventies ( see \cite{LS}). The authors consider system (\ref{NSC}) in a bounded domain $\Omega$ with
homogeneous Dirichlet
boundary conditions for $u$. Under the assumption that $u_{0}\in W^{2-\frac{2}{q},q}$ ($q>N$)
is divergence-free and vanishes
on $\p\Omega$ and that $\rho_{0}\in C^{1}(\bar{\Omega})$ is bounded away from zero, the results are the following:
\begin{itemize}
\item Global well-posedness in dimension $N=2$
\item Local well-posedness in dimension $N=3$. If in addition $u_{0}$ is small in $W^{2-\frac{2}{q},q}$, then
global well-posedness holds true.
\end{itemize}
Let us mention by passing that for the dimension $N=2$, O. Ladyzhenskaya and V. Solonnikov use a quasi-conservation
law for the $H^{1}$ norm of the velocity and get global $H^{1}$ solutions.
\\
The case of unbounded domains has been investigate by S. Itoh and A. Tani in \cite{IT}. In this framework, system (\ref{NSC})
has been shown to be locally well-posed.
In the present paper, we aim at proving similar qualitative results in the whole space $\R^{N}$ or in the torus
$\mathbb{T}^{N}$ under weaker regularity assumptions.\\
Guided in our approach by numerous works dedicated to the incompressible Navier-Stokes equation (see z.g \cite{Me}):
$$
\begin{cases}
\begin{aligned}
&\p_{t}v+v\cdot\n v-\mu\D v+\n\Pi=0,\\
&{\rm div}v=0,
\end{aligned}
\end{cases}
\leqno{(NS)}
$$
we aim at solving in the case where the data $(\rho_{0},u_{0},f)$ have critical regularity for the scaling of the equations and in particular when the initial velocity is in the same Besov spaces than Cannone, Meyer
and Planchon in \cite{CMP} for the incompressible Navier-Stokes system. It means that we want get strong solutions
when $u_{0}$ belongs to $B^{\frac{N}{p}-1}_{p,\infty}$ with $1\leq p<+\infty$.
By critical, we mean that we want to solve the system functional spaces with norm
in invariant by the changes of scales which leaves (\ref{NSC}) invariant. That approach
has been initiated by H. Fujita and T. Kato in \cite{FK}.
In the case of incompressible fluids, it is easy to see that the transformations:
$$v(t,x)\longrightarrow lv(l^{2}t,lx),\;\;\;\forall\, l\in\R$$
have that property.\\
For density-dependent incompressible fluids, one can check that the appropriate transformation are
\begin{equation}
\begin{aligned}
&(\rho_{0}(x),u_{0}(x))\longrightarrow (\rho_{0}(lx),lu_{0}(lx)),\;\;\;\forall\, l\in\R.\\
&(\rho(t,x),u(t,x),\Pi(t,x))\longrightarrow (\rho(l^{2}t,lx),lu(l^{2}t,lx),l^{2}\Pi(l^{2}t,lx)).
\end{aligned}
\label{0.5}
\end{equation}
The use of critical functional frameworks led to several new well-posedness results for incompressible fluids
(see \cite{CMP}, \cite{KT}). In the case od the density dependent incompressible fluids we want cite recent improvements by R. Danchin
in \cite{Da1}, \cite{Da2}, \cite{Da3}, \cite{Da4}, \cite{Da5}, \cite{Da6}, \cite{Da7} and of H. Abidi in \cite{Ab}
and Abidi, Paicu in \cite{AP}. All these works give results of strong solutions in finite time in critical spaces for the scaling of the equations.
More precisely R. Danchin show the existence of strong solution in finite time in \cite{Da6}, \cite{Da7}, \cite{Da8}
when the initial data check $(\rho_{0}^{-1}-1,u_{0})\in \big(B^{\N}_{2,\infty}\cap L^{\infty}\big)\times B^{\N-1}_{2,1}$
or $(\rho_{0}^{-1}-1,u_{0})\in B^{\NN}_{p,1}\times B^{\NN}_{p,1}$ with $1\leq p\leq N$. Moreover R. Danchin need of a condition of smallness on the initial data, it means that $\|\rho_{0}^{-1}-1\|_{B^{\NN}_{p,1}}$ is assumed small. More recently H. Abidi and M. Paicu
in \cite{AP} have improved this result by working with initial data in $B^{\frac{N}{p_{1}}}_{p_{1},1}\times B^{\frac{N}{p_{2}}-1}_{p_{2},1}$ with $p_{1}$ and $p_{2}$ good choosen and checking some relations. In particular they show that
we have strong solution for initial data $u_{0}$ in $B^{\frac{N}{p_{2}}-1}_{p_{2},1}$ with $1\leq p_{2}<2N$ which improves the results of R. Danchin. All these results use of crucial way the fact that the solution are Lipschitz. In particular, it explains the choice of the third index for the Besov space where $r=1$, indeed it allows a control of $\n u$ in $L^{1}_{T}(B^{\frac{N}{p}}_{p,1})$
which is embedded in $L^{1}_{T}(L^{\infty})$. This control is imperative in these works to get estimate via the transport equation on the density.\\
However the scaling of (\ref{0.5}) suggests choosing initial data $(\rho_{0},u_{0})$ in
$B^{\frac{N}{p_{1}}}_{p_{1},r^{'}}\times B^{\frac{N}{p_{2}}}_{p_{2},r}$ (see the definition of Besov spaces in the section
\ref{section2})
with $(p_{1},p_{2})\in[1,+\infty)^{2}$ and $(r,r^{'})\in[1,+\infty]^{2}$. Indeed it seems that it is not necessary by the study of the scaling of the equation to impos that $r,r^{'}=1$ as in the works of H. Abidi, R. Danchin and M. Paicu.
The goal of this article is to reach the critical case with a general third index for the Besov spaces $r$ and $r^{'}$.
More precisely in the sequel we will restrict our study to the case where the initial data $(\rho_{0},u_{0})$ and external force $f$
are such that, for some positive constant $\bar{\rho}$:
$$(\rho_{0}-\bar{\rho})\in B^{\frac{N}{p_{1}}+\e}_{p_{1},\infty}\cap L^{\infty},
\;u_{0}\in B^{\frac{N}{p_{2}}}_{p_{2},r}\;\;\mbox{and}\;\;f\in \widetilde{L}^{1}_{loc}(\R^{+},\in B^{\frac{N}{p_{2}}-1}
_{p_{2},r}).$$
with $r\in[1,+\infty]$, $\e>0$ and $p_{1}$, $p_{2}$ will verify appropriate inequalities ( for a definition of $\widetilde{L}^{1}$
see section \ref{section2}).\\
In this article we improve the result of H. Abidi, R. Danchin and M. Paicu by working with initial data in $B^{\frac{N}{p_{2}}}_{p_{2},r}$
with $r\in [1,+\infty]$. In particular we generalize to the case of the Navier-Stokes incompressible dependent density the result of Cannone-Meyer-Planchon in \cite{CMP}. For making we introduce a new idea to control the density via the transport equation when the velocity is not Lipschitz. We use some new a priori estimates on the transport equation when the velocity is only $\log$ Lipschitz. The difficulty is to deal with the loss of regularity on the density, that is why we need to ask more regularity on $\rho_{0}$ to compense this loss of regularity. The crucial point is that the density stay in a good multiplier space.
Moreover we improve the results of H. Abidi, R. Danchin and M. Paicu because we do not need to assume some condition of smallness on the initial density.
In \cite{Da1}, \cite{Da2}, \cite{Da3}, \cite{Da4}, \cite{Da5}, \cite{Da6}, \cite{Da7}, \cite{Ab}
and \cite{AP} , we need to make the additional assumption that $\rho-\bar{\rho}$ is {\it small}
in $B^{\NN}_{p,1}$.
The reason why is that we handled the elliptic operator in the momentum equation of (\ref{NSC})
as a constant coefficient second order operator plus a perturbation introduced by $\rho-\bar{\rho}$ which, if sufficiently
small, may be treated as a harmless source term. For smoother data however, the additional regularity can compensate {\it large}
perturbations of the constant coefficient operator. This fact has been used in \cite{Da4} leads to local well-posedness
results for smooth enough data with density bounded away from zero. The price to pay however is that assuming extra smoothness
precludes from using a critical framework. In fact to treat the case of large initial data on the density, we follow a idea of R. Danchin in \cite{Da} use for the case of Navier-Stokes compressible where to get good estimate on the elliptic operator in the momentum equation, we need to split the behavior high frequencies and low-middle frequencies of the viscosity terms.\\
In the present paper, we address the question of local well-posedness in the critical functional framework under the sole assumption
that the initial density is bounded away from $0$ and tends to some positive constant $\bar{\rho}$ at infinity (or
has average $\bar{\rho}$ if we consider the case of periodic boundary conditions).
To simplify the notation, we assume from now on that $\bar{\rho}=1$. Hence as long as $\rho$ does not vanish, the equations for ($a=\rho^{-1}-1$,$u$) read:
\begin{equation}
\begin{cases}
\begin{aligned}
&\p_{t}a+u\cdot\n a=0,\\
&\p_{t}u+u\cdot\n u+(1+a)(\n\Pi-\mu\D u)=f,\\
&{\rm div}u=0,\\
&(a,u)_{/t=0}=(a_{0},u_{0}).
\end{aligned}
\end{cases}
\label{0.6}
\end{equation}
One can now state our main result. The first theorem generalize the work of Cannone,Meyer, Planchon in \cite{CMP} when the third index of the Besov space for the initial data is in $[1,+\infty[$.
\begin{theorem}
Let $1\leq r<\infty$, $1\leq p_{1}<\infty$, $1<p_{2}<\infty$ and $\e>0$ such that:
$$\frac{N}{p_{1}}+\e<\frac{N}{p_{2}}+1\;\;\;\mbox{and}\;\;\;
\frac{N}{p_{2}}-1\leq \frac{N}{p_{1}}.$$
Assume that $u_{0}\in B^{\frac{N}{p_{2}}-1}_{p_{2},r}$ with ${\rm div}u_{0}=0$, $f\in \widetilde{L}^{1}_{loc}
(\R^{+},B^{\frac{N}{p_{2}}-1}_{p_{2},r})$
and $a_{0}\in B^{\frac{N}{p_{1}}+\e}_{p_{1},\infty}\cap L^{\infty}$, with
$1+a_{0}$ bounded away from zero and it exists $c>0$ such that:
$$\|a_{0}\|_{B^{\frac{N}{p_{1}}}_{p_{1},\infty}}\leq c.$$
If $\frac{1}{p_{1}}+\frac{1}{p_{2}}>\frac{1}{N}$, there exists a positive time $T$ such that system
(\ref{0.6}) has a solution
$(a,u)$ with $1+a$ bounded away from zero and:
$$
\begin{aligned}
&a\in \widetilde{C}([0,T],B^{\frac{N}{p_{1}}+\frac{\e}{2}}_{p_{1},\infty}),\;\;
u\in\big( \widetilde{C}([0,T];B^{\frac{N}{p_{2}}-1}_{p_{2},r})
\cap \widetilde{L}^{1}(0,T,B^{\frac{N}{p_{2}}+1}_{p_{2},r})\big)^{N}\\
&\hspace{7cm}\mbox{and}\;\;\n\Pi\in
\widetilde{L}^{1}(0,T,B^{\frac{N}{p_{2}}-1}_{p_{2},r}).
\end{aligned}
$$
This solution is unique when $\frac{2}{N}\leq \frac{1}{p_{1}}+\frac{1}{p_{2}}.$
\label{theo1}
\end{theorem}
\begin{remarka}
This theorem can be applied for non degenerate viscosity, but for the simplicity in the sequel the viscosity is assumed to be a constant $\mu>0$.
\end{remarka}
\begin{remarka}
As in the work of H. Abidi and M. Paicu in \cite{AP}, we are able to get strong solution when $u_{0}\in B^{\frac{N}{p_{2}}-1}_{p_{2},r} $
with $1<p_{2}\leq 2N$, it improves the result of R. Danchin in \cite{Da6} and \cite{Da7}.\\
Moreover we get weak solution with initial data very close from $(a_{0},u_{0})\in B^{1+\e}_{N,\infty}\times B^{-1}_{\infty,r}$
and $(a_{0},u_{0})\in B^{\e}_{\infty,\infty}\times B^{0}_{N,r}$. It means that in the first case we are not very far from the Koch-Tataru in \cite{KT} initial data for the velocity $u_{0}$, and in the second case we are close to just impose a condition of norm $L^{\infty}$ on $a_{0}$ which is of great interest for the system bifluid.
\end{remarka}
\begin{remarka}
In the previous theorem, we need of a condition of smallness, because when $p_{2}\ne 2$, we have extra term in our proposition \ref{linearise} which requires a condition of smallness on $a$. It seems that there is a specific structure in the scalar case when $p_{2}=2$.
\end{remarka}
In the following theorem, we improve the previous result in the specific case where $p_{2}=2$. In this case we don't need to impose condition of smallness on the initial data.
\begin{theorem}
Let $1\leq r<\infty$, $1\leq p_{1}<\infty$ and $\e>0$ such that:
$$\frac{N}{p_{1}}+\e<\frac{N}{2}+1\;\;\;\mbox{and}\;\;\;
\frac{N}{2}\leq 1+\frac{N}{p_{1}}.$$
Assume that $u_{0}\in B^{\frac{N}{2}-1}_{2,r}$ with ${\rm div}u_{0}=0$, $f\in \widetilde{L}^{1}_{loc}
(\R^{+},B^{\frac{N}{2}-1}_{2,r})$
and $a_{0}\in B^{\frac{N}{p_{1}}+\e}_{p_{1},\infty}\cap L^{\infty}$, with
$1+a_{0}$ bounded away from zero.
There exists a positive time $T$ such that system
(\ref{0.6}) has a solution
$(a,u)$ with $1+a$ bounded away from zero and:
$$
\begin{aligned}
&a\in \widetilde{C}([0,T],B^{\frac{N}{p_{1}}+\frac{\e}{2}}_{p_{1},\infty}),\;\;
u\in\big( \widetilde{C}([0,T];B^{\frac{N}{2}-1}_{2,r})
\cap \widetilde{L}^{1}(0,T,B^{\frac{N}{2}+1}_{2,r})\big)^{N}\\
&\hspace{7cm}\mbox{and}\;\;\n\Pi\in
\widetilde{L}^{1}(0,T,B^{\frac{N}{2}-1}_{2,r}).
\end{aligned}
$$
This solution is unique when $\frac{2}{N}\leq \frac{1}{p_{1}}+\frac{1}{2}.$
\label{theo11}
\end{theorem}
In the following theorem we choose $r=+\infty$, we have to treat the case of a linear loss of regularity
on the density $\rho$ which depends of the solution $u$. This fact come from your estimate on the transport equation when $\n u\in\widetilde{L}_{T}^{1}(B^{\frac{N}{p_{2}}}_{p_{2},\infty})$, in this case we have a loss of regularity on the density which depends of
this quantity.
\begin{theorem}
Let $1\leq p_{1}<\infty$, $1<p_{2}<\infty$,  and $\e>0$ such that:
$$\frac{N}{p_{1}}+\e<\frac{N}{p_{2}}+1\;\;\;\mbox{and}\;\;\;
\frac{N}{p_{2}}\leq 1+\frac{N}{p_{1}}.$$
Assume that $u_{0}\in B^{\frac{N}{p_{2}}-1}_{p_{2},\infty}$ with ${\rm div}u_{0}=0$, $f\in \widetilde{L}^{1}_{loc}
(\R^{+},B^{\frac{N}{p_{2}}-1}_{p_{2},\infty})$
and $a_{0}\in B^{\frac{N}{p_{1}}+\e}_{p_{1},\infty}\cap L^{\infty}$, with
$1+a_{0}$ bounded away from zero and it exists $c>0$ such that:
$$\|a_{0}\|_{B^{\frac{N}{p_{1}}}_{p_{1},\infty}}\leq c.$$
If $\frac{1}{p_{1}}+\frac{1}{p_{2}}>\frac{1}{N}$, there exists a positive time $T$ such that system
(\ref{0.6}) has a solution
$(a,u)$ with $1+a$ bounded away from zero and:
$$
\begin{aligned}
&a\in \widetilde{C}([0,T],B^{\sigma(T)}_{p_{1},\infty}),\;\;
u\in\big( \widetilde{C}([0,T];B^{\frac{N}{p_{2}}-1}_{p_{2},\infty}
\cap \widetilde{L}^{1}(0,T,B^{\frac{N}{p_{2}}+1}_{p_{2},\infty})\big)^{N}\\
&\hspace{7cm}\mbox{and}\;\;\n\Pi\in\widetilde{L}^{1}(0,T,B^{\frac{N}{p_{2}}-1}_{p_{2},\infty}),
.\end{aligned}
$$
with:
$$\sigma(T)=\frac{N}{p_{1}}+\e-\lambda\|u\|_{\widetilde{L}^{1}_{T}(B^{\frac{N}{p_{2}}+1}_{p_{2},\infty})}
+\lambda\int^{T}_{0}W(t^{'})dt^{'},$$
for any $\lambda>0$ and any nonegative integrable function $W$ over $[0,T]$ such that $\sigma(T)\geq-1-N\min
(\frac{1}{p_{1}},\frac{1}{p_{1}^{'}})$. 
This solution is unique when $\frac{2}{N}\leq \frac{1}{p_{1}}+\frac{1}{p_{2}}.$
\label{theo2}
\end{theorem}
In the following theorem, we generalize the previous result with large initial data for the density when $p_{2}=2$.
\begin{theorem}
Let $1\leq p_{1}<\infty$  and $\e>0$ such that:
$$\frac{N}{p_{1}}+\e<\frac{N}{2}+1\;\;\;\mbox{and}\;\;\;
\frac{N}{2}\leq 1+\frac{N}{p_{1}}.$$
Assume that $u_{0}\in B^{\frac{N}{2}-1}_{2,\infty}$ with ${\rm div}u_{0}=0$, $f\in \widetilde{L}^{1}_{loc}
(\R^{+},B^{\frac{N}{2}-1}_{2,\infty})$
and $a_{0}\in B^{\frac{N}{p_{1}}+\e}_{p_{1},\infty}\cap L^{\infty}$, with
$1+a_{0}$ bounded away from zero.
There exists a positive time $T$ such that system
(\ref{0.6}) has a solution
$(a,u)$ with $1+a$ bounded away from zero and:
$$
\begin{aligned}
&a\in \widetilde{C}([0,T],B^{\sigma(T)}_{p_{1},\infty}),\;\;
u\in\big( \widetilde{C}([0,T];B^{\frac{N}{2}-1}_{2,\infty}
\cap \widetilde{L}^{1}(0,T,B^{\frac{N}{2}+1}_{2,\infty})\big)^{N}\\
&\hspace{7cm}\mbox{and}\;\;\n\Pi\in\widetilde{L}^{1}(0,T,B^{\frac{N}{2}-1}_{2,\infty}),
.\end{aligned}
$$
with:
$$\sigma(T)=\frac{N}{p_{1}}+\e-\lambda\|u\|_{\widetilde{L}^{1}_{T}(B^{\frac{N}{2}+1}_{2,\infty})}
+\lambda\int^{T}_{0}W(t^{'})dt^{'},$$
for any $\lambda>0$ and any nonegative integrable function $W$ over $[0,T]$ such that $\sigma(T)\geq-1-N\min
(\frac{1}{p_{1}},\frac{1}{p_{1}^{'}})$. 
This solution is unique when $\frac{2}{N}\leq \frac{1}{p_{1}}+\frac{1}{2}.$
\label{theo22}
\end{theorem}
The key of the theorems \ref{theo1}, \ref{theo11}, \ref{theo2} and \ref{theo22} is new estimate for transport equation for $u$ which is not considered Lipschitz. In this case we have to pay a loss of regularity on the density.
$\rho$.
The basic idea is to deal with this loss of regularity by conserving $a$ in $\widetilde{C}_{T}(B^{\frac{N}{p_{1}}}_{p_{1},+\infty})\cap L^{\infty}$ which have good properties of multiplier.
Moreover to avoid as in the work s of H. Abidi, R. Danchin and M. Paicu conditions of smallness on the initial density, the basic idea is that having coefficients in $a$ in $\widetilde{C}_{T}(B^{\frac{N}{p_{1}}}_{p_{1},+\infty})$ provides us with some uniform
decay on the high frequencies of the variable coefficients so that the elliptic operator may be
considered as an operator with
{\it smooth} coefficients (of the type which has been investigated in \cite{Da5}) plus a small error term.\\
Our paper is structured as follows. In the section \ref{section2}, we give a few notation and briefly introduce the basic Fourier
analysis
techniques needed to prove our result. Section \ref{section3} and \ref{section4} are devoted to the proof of key estimates for the linearized
system (\ref{0.6}) in particular the elliptic operator of the momentum equation with variable coefficients and the transport equation when the velocity is not assumed Lipschitz. In section \ref{section5}, we prove the existenceof solution for theorem \ref{theo1} and \ref{theo11} whereas section
\ref{section6}
is devoted to the proof of uniqueness. Finally in section \ref{section7}, we briefly show how to prove theorem  \ref{theo2} and \ref{theo22}.
Elliptic and technical estimates commutator
are postponed in an appendix.
\section{Littlewood-Paley theory and Besov spaces}
\label{section2}
Throughout the paper, $C$ stands for a constant whose exact meaning depends on the context. The notation $A\lesssim B$ means
that $A\leq CB$.
For all Banach space $X$, we denote by $C([0,T],X)$ the set of continuous functions on $[0,T]$ with values in $X$.
For $p\in[1,+\infty]$, the notation $L^{p}(0,T,X)$ or $L^{p}_{T}(X)$ stands for the set of measurable functions on $(0,T)$
with values in $X$ such that $t\rightarrow\|f(t)\|_{X}$ belongs to $L^{p}(0,T)$.
\subsection{Littlewood-Paley decomposition}
Littlewood-Paley decomposition  corresponds to a dyadic
decomposition  of the space in Fourier variables.
Let $\alpha>1$ and $(\va,\chi)$ be a couple of smooth functions valued in $[0,1]$, such that $\va$ is supported in the shell
supported in
$\{\xi\in\R^{N}/\alpha^{-1}\leq|\xi|\leq2\alpha\}$, $\chi$ is supported in the ball $\{\xi\in\R^{N}/|\xi|\leq\alpha\}$
such that:
$$\forall\xi\in\R^{N},\;\;\;\chi(\xi)+\sum_{l\in\mathbb{N}}\varphi(2^{-l}\xi)=1.$$
Denoting $h={\cal{F}}^{-1}\varphi$, we then define the dyadic
blocks by:
$$
\begin{aligned}
&\D_{l}u=0\;\;\;\mbox{if}\;\;l\leq-2,\\
&\D_{-1}u=\chi(D)u=\widetilde{h}*u\;\;\;\mbox{with}\;\;\widetilde{h}={\cal F}^{-1}\chi,\\
&\D_{l}u=\varphi(2^{-l}D)u=2^{lN}\int_{\R^{N}}h(2^{l}y)u(x-y)dy\;\;\;\mbox{with}\;\;h={\cal F}^{-1}\chi,\;\;\mbox{if}\;\;l\geq0,\\
&S_{l}u=\sum_{k\leq
l-1}\D_{k}u\,.
\end{aligned}
$$
Formally, one can write that:
$u=\sum_{k\in\mathbb{Z}}\D_{k}u\,.$
This decomposition is called nonhomogeneous Littlewood-Paley
decomposition. 
\subsection{Nonhomogeneous Besov spaces and first properties}
\begin{definition}
For
$s\in\R,\,\,p\in[1,+\infty],\,\,q\in[1,+\infty],\,\,\mbox{and}\,\,u\in{\cal{S}}^{'}(\R^{N})$
we set:
$$\|u\|_{B^{s}_{p,q}}=(\sum_{l\in\mathbb{Z}}(2^{ls}\|\D_{l}u\|_{L^{p}})^{q})^{\frac{1}{q}}.$$
The Besov space $B^{s}_{p,q}$ is the set of temperate distribution $u$ such that $\|u\|_{B^{s}_{p,q}}<+\infty$.
\end{definition}
\begin{remarka}The above definition is a natural generalization of the
nonhomogeneous Sobolev and H$\ddot{\mbox{o}}$lder spaces: one can show
that $B^{s}_{\infty,\infty}$ is the nonhomogeneous
H$\ddot{\mbox{o}}$lder space $C^{s}$ and that $B^{s}_{2,2}$ is
the nonhomogeneous space $H^{s}$.
\end{remarka}
\begin{proposition}
\label{derivation,interpolation}
The following properties holds:
\begin{enumerate}
\item there exists a constant universal $C$
such that:\\
$C^{-1}\|u\|_{B^{s}_{p,r}}\leq\|\n u\|_{B^{s-1}_{p,r}}\leq
C\|u\|_{B^{s}_{p,r}}.$
\item If
$p_{1}<p_{2}$ and $r_{1}\leq r_{2}$ then $B^{s}_{p_{1},r_{1}}\hookrightarrow
B^{s-N(1/p_{1}-1/p_{2})}_{p_{2},r_{2}}$.
\item
 $(B^{s_{1}}_{p,r},B^{s_{2}}_{p,r})_{\theta,r^{'}}=B^{\theta
s_{1}+(1-\theta)s_{2}}_{p,r^{'}}$.
\end{enumerate}
\end{proposition}
Let now recall a few product laws in Besov spaces coming directly from the paradifferential calculus of J-M. Bony
(see \cite{Bo}) and rewrite on a generalized form in \cite{AP} by H. Abidi and M. Paicu (in this article the results are written
in the case of homogeneous sapces but it can easily generalize for the nonhomogeneous Besov spaces).
\begin{proposition}
\label{produit1}
We have the following laws of product:
\begin{itemize}
\item For all $s\in\R$, $(p,r)\in[1,+\infty]^{2}$ we have:
\begin{equation}
\|uv\|_{\widetilde{B}^{s}_{p,r}}\leq
C(\|u\|_{L^{\infty}}\|v\|_{B^{s}_{p,r}}+\|v\|_{L^{\infty}}\|u\|_{B^{s}_{p,r}})\,.
\label{2.2}
\end{equation}
\item Let $(p,p_{1},p_{2},r,\lambda_{1},\lambda_{2})\in[1,+\infty]^{2}$ such that:$\frac{1}{p}\leq\frac{1}{p_{1}}+\frac{1}{p_{2}}$,
$p_{1}\leq\lambda_{2}$, $p_{2}\leq\lambda_{1}$, $\frac{1}{p}\leq\frac{1}{p_{1}}+\frac{1}{\lambda_{1}}$ and
$\frac{1}{p}\leq\frac{1}{p_{2}}+\frac{1}{\lambda_{2}}$. We have then the following inequalities:\\
if $s_{1}+s_{2}+N\inf(0,1-\frac{1}{p_{1}}-\frac{1}{p_{2}})>0$, $s_{1}+\frac{N}{\lambda_{2}}<\frac{N}{p_{1}}$ and
$s_{2}+\frac{N}{\lambda_{1}}<\frac{N}{p_{2}}$ then:
\begin{equation}
\|uv\|_{B^{s_{1}+s_{2}-N(\frac{1}{p_{1}}+\frac{1}{p_{2}}-\frac{1}{p})}_{p,r}}\lesssim\|u\|_{B^{s_{1}}_{p_{1},r}}
\|v\|_{B^{s_{2}}_{p_{2},\infty}},
\label{2.3}
\end{equation}
when $s_{1}+\frac{N}{\lambda_{2}}=\frac{N}{p_{1}}$ (resp $s_{2}+\frac{N}{\lambda_{1}}=\frac{N}{p_{2}}$) we replace
$\|u\|_{B^{s_{1}}_{p_{1},r}}\|v\|_{B^{s_{2}}_{p_{2},\infty}}$ (resp $\|v\|_{B^{s_{2}}_{p_{2},\infty}}$) by
$\|u\|_{B^{s_{1}}_{p_{1},1}}\|v\|_{B^{s_{2}}_{p_{2},r}}$ (resp $\|v\|_{B^{s_{2}}_{p_{2},\infty}\cap L^{\infty}}$),
if $s_{1}+\frac{N}{\lambda_{2}}=\frac{N}{p_{1}}$ and $s_{2}+\frac{N}{\lambda_{1}}=\frac{N}{p_{2}}$ we take $r=1$.
\\
If $s_{1}+s_{2}=0$, $s_{1}\in(\frac{N}{\lambda_{1}}-\frac{N}{p_{2}},\frac{N}{p_{1}}-\frac{N}{\lambda_{2}}]$ and
$\frac{1}{p_{1}}+\frac{1}{p_{2}}\leq 1$ then:
\begin{equation}
\|uv\|_{B^{-N(\frac{1}{p_{1}}+\frac{1}{p_{2}}-\frac{1}{p})}_{p,\infty}}\lesssim\|u\|_{B^{s_{1}}_{p_{1},1}}
\|v\|_{B^{s_{2}}_{p_{2},\infty}}.
\label{2.4}
\end{equation}
If $|s|<\NN$ for $p\geq2$ and $-\frac{N}{p^{'}}<s<\NN$ else, we have:
\begin{equation}
\|uv\|_{B^{s}_{p,r}}\leq C\|u\|_{B^{s}_{p,r}}\|v\|_{B^{\NN}_{p,\infty}\cap L^{\infty}}.
\label{2.5}
\end{equation}
\end{itemize}
\end{proposition}
\begin{remarka}
In the sequel $p$ will be either $p_{1}$ or $p_{2}$ and in this case $\frac{1}{\lambda}=\frac{1}{p_{1}}-\frac{1}{p_{2}}$
if $p_{1}\leq p_{2}$, resp $\frac{1}{\lambda}=\frac{1}{p_{2}}-\frac{1}{p_{1}}$
if $p_{2}\leq p_{1}$.
\end{remarka}
\begin{corollaire}
\label{produit2}
Let $r\in [1,+\infty]$, $1\leq p\leq p_{1}\leq +\infty$ and $s$ such that:
\begin{itemize}
\item $s\in(-\frac{N}{p_{1}},\frac{N}{p_{1}})$ if $\frac{1}{p}+\frac{1}{p_{1}}\leq 1$,
\item $s\in(-\frac{N}{p_{1}}+N(\frac{1}{p}+\frac{1}{p_{1}}-1),\frac{N}{p_{1}})$ if $\frac{1}{p}+\frac{1}{p_{1}}> 1$,
\end{itemize}
then we have if $u\in B^{s}_{p,r}$ and $v\in B^{\frac{N}{p_{1}}}_{p_{1},\infty}\cap L^{\infty}$:
$$\|uv\|_{B^{s}_{p,r}}\leq C\|u\|_{B^{s}_{p,r}}\|v\|_{B^{\frac{N}{p_{1}}}_{p_{1},\infty}\cap L^{\infty}}.$$
\end{corollaire}
The study of non stationary PDE's requires space of type $L^{\rho}(0,T,X)$ for appropriate Banach spaces $X$. In our case, we
expect $X$ to be a Besov space, so that it is natural to localize the equation through Littlewood-Payley decomposition. But, in doing so, we obtain
bounds in spaces which are not type $L^{\rho}(0,T,X)$ (except if $r=p$).
We are now going to
define the spaces of Chemin-Lerner in which we will work, which are
a refinement of the spaces
$L_{T}^{\rho}(B^{s}_{p,r})$.
$\hspace{15cm}$
\begin{definition}
Let $\rho\in[1,+\infty]$, $T\in[1,+\infty]$ and $s_{1}\in\R$. We set:
$$\|u\|_{\widetilde{L}^{\rho}_{T}(B^{s_{1}}_{p,r})}=
\big(\sum_{l\in\mathbb{Z}}2^{lrs_{1}}\|\D_{l}u(t)\|_{L^{\rho}(L^{p})}^{r}\big)^{\frac{1}{r}}\,.$$
We then define the space $\widetilde{L}^{\rho}_{T}(B^{s_{1}}_{p,r})$ as the set of temperate distribution $u$ over
$(0,T)\times\R^{N}$ such that 
$\|u\|_{\widetilde{L}^{\rho}_{T}(B^{s_{1}}_{p,r})}<+\infty$.
\end{definition}
We set $\widetilde{C}_{T}(\widetilde{B}^{s_{1}}_{p,r})=\widetilde{L}^{\infty}_{T}(\widetilde{B}^{s_{1}}_{p,r})\cap
{\cal C}([0,T],B^{s_{1}}_{p,r})$.
Let us emphasize that, according to Minkowski inequality, we have:
$$\|u\|_{\widetilde{L}^{\rho}_{T}(B^{s_{1}}_{p,r})}\leq\|u\|_{L^{\rho}_{T}(B^{s_{1}}_{p,r})}\;\;\mbox{if}\;\;r\geq\rho
,\;\;\;\|u\|_{\widetilde{L}^{\rho}_{T}(B^{s_{1}}_{p,r})}\geq\|u\|_{L^{\rho}_{T}(B^{s_{1}}_{p,r})}\;\;\mbox{if}\;\;r\leq\rho
.$$
\begin{remarka}
It is easy to generalize proposition \ref{produit1},
to $\widetilde{L}^{\rho}_{T}(B^{s_{1}}_{p,r})$ spaces. The indices $s_{1}$, $p$, $r$
behave just as in the stationary case whereas the time exponent $\rho$ behaves according to H\"older inequality.
\end{remarka}
Here we recall a result of interpolation which explains the link
of the space $B^{s}_{p,1}$ with the space $B^{s}_{p,\infty}$, see
\cite{Da1}.
\begin{proposition}
\label{interpolationlog}
There exists a constant $C$ such that for all $s\in\R$, $\e>0$ and
$1\leq p<+\infty$,
$$\|u\|_{\widetilde{L}_{T}^{\rho}(B^{s}_{p,1})}\leq C\frac{1+\e}{\e}\|u\|_{\widetilde{L}_{T}^{\rho}(B^{s}_{p,\infty})}
\biggl(1+\log\frac{\|u\|_{\widetilde{L}_{T}^{\rho}(B^{s+\e}_{p,\infty})}}
{\|u\|_{\widetilde{L}_{T}^{\rho}(B^{s}_{p,\infty})}}\biggl).$$ \label{5Yudov}
\end{proposition}
\begin{definition}
\label{defChemin1}
Let $\Gamma$ be an increasing function on $[1,+\infty[$. We denote by $B_{\Gamma}(\R^{N})$ the set of bounded
real valued functions $u$ over $\R^{N}$ such that:
$$\|u\|_{B_{\Gamma}}=\|u\|_{L^{\infty}}+\sup_{j\geq0}\frac{\|\n S_{j}u\|_{L^{\infty}}}{\Gamma(2^{j})}<+\infty.$$
\end{definition}
We give here a proposition concerning these spaces showed by J-Y. Chemin, see \cite{BCD}.
\begin{proposition}
\label{propChemin2}
Let $\e>0$ and $u\in\widetilde{L}_{T}^{1}(B^{\NN+1}_{p,r})$ then we have $u\in L^{1}_{T}(B_{\Gamma}(\R^{N}))$ with
$\Gamma(t)=(-\log t)^{1+\e-\frac{1}{r}}$ for $0\leq t\leq 1$ .
\end{proposition}
\section{Estimates for parabolic system with variable coefficients}
\label{section3}
In this section, the following linearization of the momentum equation is studied:
\begin{equation}
\begin{cases}
\begin{aligned}
&\p_{t}u+b(\n\Pi-\mu\D u)=f+g,\\
&{\rm div}u=0,\\
&u_{/t=0}=u_{0}.
\end{aligned}
\end{cases}
\label{5}
\end{equation}
where $b$, $f$, $g$ and $u_{0}$ are given.
Above $u$ is the unknown function.
We assume that $u_{0}\in B^{s}_{p,r}$ and $f\in \widetilde{L}^{1}(0,T;B^{s}_{p,r})$,
that $b$ is bounded by below by a positive constant $\underline{b}$ and that
$a=b-1$ belongs to $\widetilde{L}^{\infty}(0,T;B^{\frac{N}{p_{1}}+\alpha}_{p_{1},\infty})\cap L^{\infty}$.
In the present subsection, we aim at proving a priori estimates for (\ref{5}) in the framework of nonhomogeneous
Besov spaces. Before stating our results let us introduce the following notation:
\begin{equation}
\begin{aligned}
{\cal A}_{T}=&1+\underline{b}^{-1}\|\n b\|_{\widetilde{L}^{\infty}(B^{\frac{N}{p_{1}}+\alpha-1}_{p_{1},\infty})}\;\;\;\mbox{with}\;\;\alpha>0.
\end{aligned}
\label{notation16}
\end{equation}
\begin{proposition}
Let $\underline{\nu}=\underline{b}\mu$ and $(p,p_{1})\in [1,+\infty]$.
\begin{itemize}
\item If $p_{1}>p$ we assume that
$s\in(-\frac{N}{p_{1}},\frac{N}{p_{1}})$ if $\frac{1}{p}+\frac{1}{p_{1}}\leq1$ and
$s\in(-\frac{N}{p_{1}}+N(\frac{1}{p}+\frac{1}{p_{1}}-1),\frac{N}{p_{1}})$ if $\frac{1}{p}+\frac{1}{p_{1}}>1$.
\item If $p_{1}\leq p$ then we suppose that
$s\in(-\frac{N}{p},\frac{N}{p})$ if $p\geq2$ and
$s\in(-\frac{N}{p^{'}},\frac{N}{p})$ if $p<2$.
\end{itemize}
If $p\ne 2$ we need to assume than there exists $c>0$ such that:
$$\|a\|_{\widetilde{L}^{\infty}(B^{\frac{N}{p_{1}}+\alpha-1}_{p_{1},\infty})}\leq c.$$
Let $m\in\mathbb{Z}$
be such that $b_{m}=1+S_{m}a$ satisfies:
\begin{equation}
\inf_{(t,x)\in[0,T)\times\R^{N}}b_{m}(t,x)\geq\frac{\underline{b}}{2}.
\label{6}
\end{equation}
There exist three constants $c$, $C$ and $\kappa$ (with $c$, $C$, depending only on $N$ and on $s$, and $\kappa$ universal) such that if in addition we have:
\begin{equation}
\|a-S_{m}a\|_{\widetilde{L}^{\infty}(0,T;B^{\frac{N}{p_{1}}}_{p_{1},\infty})\cap L^{\infty}}\leq c\frac{\underline{\nu}}{\mu}
\label{7}
\end{equation}
then setting:
$$Z_{m}(t)=2^{2m\alpha}\mu^{2}\underline{\nu}^{-1}\int^{t}_{0}\|a\|^{2}_{B^{\frac{N}{p_{1}}}_{p_{1},\infty}\cap L^{\infty}}d\tau,$$
Let $\alpha^{'}>0$ checking $\alpha^{'}\leq\min(1,\alpha,\frac{s-2+\frac{2}{m}}{2})$.
We have for all $t\in[0,T]$ and $\kappa=\frac{s}{\alpha^{'}}$:
\begin{equation}
\begin{aligned}
&\|u\|_{\widetilde{L}^{\infty}_{T}(B^{s}_{p,r})}+\kappa\underline{\nu}
\|u\|_{\widetilde{L}^{1}_{T}(B^{s+2}_{p,r})}\leq e^{CZ_{m}(T)}\biggl(\|u_{0}\|_{B^{s}_{p,r}}+{\cal A}_{T}^{\kappa}(
\|{\cal P}f\|_{\widetilde{L}^{1}_{T}(B^{s}_{p,r})}\\
&\hspace{2,5cm}+\underline{\mu}^{\frac{1}{m}}\|{\cal P}g\|_{\widetilde{L}^{m}_{T}(B^{s-2+\frac{2}{m}}_{p,r})}+\underline{\mu}^{\frac{1}{m}}(\frac{\underline{\nu}(p-1)}{p}){\cal A}_{T}\|u\|_{\widetilde{L}^{1}_{T}
(B^{s+2-\alpha^{'}}_{p,r})}\big)\biggl).
\end{aligned}
\label{b16}
\end{equation}
Moreover we have $\n\Pi=\n\Pi_{1}+\n\Pi_{2}$ with:
\begin{equation}
\begin{aligned}
&\underline{b}\|\n \Pi_{1}\|_{\widetilde{L}^{1}_{T}(B^{s}_{p,r})}\leq {\cal A}^{\kappa}_{T}\|{\cal P}f\|_{\widetilde{L}^{1}_{T}(B^{s}_{p,r})},\\
&\underline{b}\|\n \Pi_{2}\|_{\widetilde{L}^{1}_{T}(B^{s}_{p,r})}\leq {\cal A}^{\kappa}_{T}\big(\|{\cal Q}g\|_{\widetilde{L}^{m}_{T}(B^{s-2+\frac{2}{m}}_{p,r})}+\mu\|a\|_{\widetilde{L}^{\infty}_{T}(B^{\frac{N}{p_{1}}+\alpha}_{p_{1},+\infty})}
\|\D u\|_{\widetilde{L}^{m}_{T}(B^{s-2+\frac{2}{m}}_{p,r})}\big).
\end{aligned}
\label{b17}
\end{equation}
\label{linearise}
\end{proposition}
\begin{remarka}
Let us stress the fact that if $a\in \widetilde{L}^{\infty}((0,T)\times B^{\frac{N}{p_{1}}}_{p_{1},\infty})$ then assumption (\ref{6}) and
(\ref{7}) are satisfied for $m$ large enough. This will be used in the proof of theorem \ref{theo11} and \ref{theo22}.
Indeed, according to Bernstein inequality
for $m$ large enough 9\ref{6}) and (\ref{7}) are satisfied.
\end{remarka}
Proving proposition \ref{linearise} in the case $b=cste$ is not too involved as one
can easily get rid of the pressure by taking advantage of the Leray projector ${\cal P}$ on solenoidal vector-fields.
Then system (\ref{5}) reduce to a linear $\psi$DO which may be easily solved by mean of energy estimates.
In our case where $b$ is not assumed to be a constant, getting rid of the pressure will still be an appropriate strategy.
This may be achieved by applying the operator ${\rm div}$ to (\ref{5}). Indeed by doing so, we see that the pressure
solves the elliptic equation:
\begin{equation}
{\rm div}(b\n\Pi)={\rm div}F.
\label{ellipticequation}
\end{equation}
with $F=f+g+\mu a\D u.$
Therefore denoting by ${\cal H}_{b}$ the linear operator $F\rightarrow\n\Pi$, system (\ref{5}) reduces to a linear ODE in
Banach spaces.
Actually, due to the consideration of two forcing terms $f$ and $g$ with different regularities, the pressure
has to be split into two parts, namely $\Pi=\Pi_{1}+\Pi_{2}$ with:
\begin{equation}
{\rm div}(b\Pi_{1})={\rm div}f
\label{191}
\end{equation}
\begin{equation}
{\rm div}(b\Pi_{2})={\rm div}H\;\;\mbox{and}\;\;H=g+\mu a\D u.
\label{201}
\end{equation}
{\bf Proof of proposition \ref{linearise}:}\\
\\
Let us first rewrite (\ref{5}) as follows:
\begin{equation}
\begin{cases}
\begin{aligned}
&\p_{t}u-b_{m}\mu\D u+b\n\Pi=f+g+E_{m},\\
&{\rm div}u=0,\\
&u_{t=0}=u_{0}.
\end{aligned}
\end{cases}
\label{8}
\end{equation}
with $E_{m}=\mu\D u\,(\mbox{Id}-S_{m})a$ and $b_{m}=1+S_{m}a$.
Note that by using corollary \ref{produit2} and as $-\frac{N}{p_{1}}<s<\frac{N}{p_{1}}$ for $p\geq2$ or $\frac{N}{p_{1}^{'}}<s<\frac{N}{p_{1}}$ else, the error term $E_{m}$ may be estimated by:
\begin{equation}
\|E_{m}\|_{B^{s}_{p,r}}\lesssim\|a-S_{m}a\|_{B^{\frac{N}{p_{1}}}_{p_{1},\infty}\cap L^{\infty}}\|D^{2}u\|_{B^{s}_{p,r}}.
\label{9}
\end{equation}
Now applying operator $\D_{q}$ and next operator of free divergence yield ${\cal P}$ to momentum equation (\ref{8}) yields:
\begin{equation}
\begin{aligned}
\frac{d}{dt}u_{q}-\mu{\rm div}(b_{m}\n u_{q})={\cal P}f_{q}+{\cal P}g_{q}+\D_{q}{\cal P}E_{m}+\widetilde{R}_{q}-\D_{q}{\cal P}(a\n\Pi),
\end{aligned}
\end{equation}
where we denote by $u_{q}=\D_{q}u$ and
with:
$$
\begin{aligned}
&\widetilde{R}_{q}=\widetilde{R}^{1}_{q}+\widetilde{R}^{2}_{q}.
\end{aligned}
$$
where:
$$
\begin{aligned}
&\widetilde{R}^{1}_{q}=\mu\big({\cal P}\D_{q}(b_{m}\D u)-{\cal P}{\rm div}(b_{m}\n u_{q})\big),\\
&\widetilde{R}^{2}_{q}=\mu\big({\cal P}{\rm div}(b_{m}\n u_{q})-{\rm div}(b_{m}\n u_{q})\big)=-\mu {\cal Q}{\rm div}(S_{m}a\n u_{q}).
\end{aligned}
$$
where ${\cal Q}$ is the gradient yield projector.
\subsubsection*{case $p\ne 2$}
Next multiplying both sides by $|u_{q}|^{p-2}u_{q}$, and integrating by parts in the second, third and last term in the
left-hand side, we get by using Bony decomposition (for the notation see \cite{Da1}):
\begin{equation}
\begin{aligned}
&\frac{1}{p}\frac{d}{dt}\|u_{q}\|_{L^{p}}^{p}+\mu\int_{\R^{N}}
b_{m}|\n u_{q}|^{2}|u_{q}|^{p-2}dx+\mu\int_{\R^{N}}b_{m}|u_{q}|^{p-4}|\n|u|^{2}|^{2}dx\\
&\leq\|u_{q}\|^{p-1}_{L^{p}}(\|{\cal P}f_{q}\|_{L^{p}}+\|{\cal P}g_{q}\|_{L^{p}}
+\|\widetilde{R}_{q}\|_{L^{p}}+\|\D_{q}(T_{\n a}\Pi)\|_{L^{p}}+2^{q}\|\D_{q}(T_{a}\Pi)\|_{L^{p}}\\
&\hspace{7,5cm}+\|\D_{q}(T^{'}_{\n \Pi}a)\|_{L^{p}}+\|{\cal P}\D_{q}E_{m}\|_{L^{p}}).
\end{aligned}
\label{energieDanch}
\end{equation}
Indeed we have as ${\rm div}u=0$ and by using Bony's decomposition and by performing an integration by parts:
$$
\begin{aligned}
&\int_{\R^{N}}\D_{q}(a\n\Pi)|u_{q}|^{p-2}u_{q}=\int_{\R^{N}}\D_{q}(T^{'}_{\n \Pi}a)|u_{q}|^{p-2}u_{q}dx-
\int_{\R^{N}}\D_{q}(T_{\n a}\Pi)|u_{q}|^{p-2}u_{q}dx\\
&\hspace{8,7cm}-\int_{\R^{N}}\D_{q}(T_{a}\Pi){\rm div}(|u_{q}|^{p-2} u_{q})dx.
\end{aligned}
$$
Next we have:
$$\n(|u_{q}|^{p-2})\cdot u_{q}=(p-2)|u_{q}|^{p-4}\sum_{i,k}u_{q}^{k}\p_{i}u_{q}^{k}u_{q}^{i},$$
and by H\"older's and Berstein's inequalities:
$$\|\n(|u_{q}|^{p-2})\cdot u_{q}\|_{L^{\frac{p}{p-1}}}\leq C(p-2)2^{q}\|u_{q}\|_{L^{p}}^{p-1}.$$
Next from inequality (\ref{energieDanch}), we get by using lemma A5 in \cite{Da6}:
$$
\begin{aligned}
&\frac{1}{p}\frac{d}{dt}\|u_{q}\|_{L^{p}}^{p}+\frac{\underline{\nu}(p-1)}{p^{2}}2^{2q}\|u_{q}\|_{L^{p}}^{p}\leq
\|u_{q}\|^{p-1}_{L^{p}}\big(\|{\cal P}f_{q}\|_{L^{p}}+\|{\cal P}g_{q}\|_{L^{p}}+\|{\cal P}\D_{q}E_{m}\|_{L^{p}}
\\
&\hspace{2,5cm}+\|\D_{q}(T_{\n a}\Pi)\|_{L^{p}}+2^{q}\|\D_{q}(T_{a}\Pi)\|_{L^{p}}+\|\D_{q}(T^{'}_{\n \Pi}a)\|_{L^{p}}+\|\widetilde{R}_{q}\|_{L^{p}}\big),
\end{aligned}
$$
Therefore, elementary computation yield (at least formally):
$$
\begin{aligned}
&e^{-\frac{\underline{\nu}(p-1)}{p^{2}}2^{2q}t}\frac{d}{dt}
\big(e^{\frac{\underline{\nu}(p-1)}{p^{2}}2^{2q}t}\|u_{q}\|_{L^{p}}\big)\lesssim
\|{\cal P}f_{q}\|_{L^{p}}+\|{\cal P}g_{q}\|_{L^{p}}+\|{\cal P}\D_{q}E_{m}\|_{L^{p}}\\
&\hspace{2cm}+\|\D_{q}(T_{\n a}\Pi)\|_{L^{p}}+2^{q}\|\D_{q}(T_{a}\Pi)\|_{L^{p}}+\|\D_{q}(T^{'}_{\n \Pi}a)\|_{L^{p}}+
\|\widetilde{R}_{q}\|_{L^{p}}.
\end{aligned}
$$
We thus have:
$$
\begin{aligned}
&\|u_{q}(t)\|_{L^{p}}\lesssim e^{-\frac{\underline{\nu}(p-1)}{p^{2}}2^{2q}t}\|\D_{q}u_{0}\|_{L^{p}}+
\int^{t}_{0}e^{-\frac{\underline{\nu}(p-1)}{p^{2}}2^{2q}(t-\tau)}\big(\|{\cal P}f_{q}\|_{L^{p}} +\|{\cal P}g_{q}\|_{L^{p}}+\\
&\|{\cal P}\D_{q}E_{m}\|_{L^{p}}+\|\D_{q}(T_{\n a}\Pi)\|_{L^{p}}+2^{q}\|\D_{q}(T_{a}\Pi)\|_{L^{p}}+\|\D_{q}(T^{'}_{\n \Pi}a)\|_{L^{p}}+
\|\widetilde{R}_{q}\|_{L^{p}}\big)(\tau)d\tau,
\end{aligned}
$$
which leads for all $q\geq -1$, after performing a time integration and using convolution inequalities to:
\begin{equation}
\begin{aligned}
&(\frac{\underline{\nu}(p-1)}{p^{2}})^{\frac{1}{m}}2^{\frac{2q}{m}}\|u_{q}\|_{L^{m}_{T}(L^{p})}\lesssim
\|\D_{q}u_{0}\|_{L^{p}}+\|{\cal P}f_{q}\|_{L^{1}_{T}(L^{p})}+\|\D_{q}(T_{\n a}\Pi_{1})\|_{L^{1}_{T}(L^{p})}\\
&+2^{q}\|\D_{q}(T_{a}\Pi_{1})\|_{L^{1}_{T}(L^{p})}+\|\D_{q}(T^{'}_{\n \Pi_{1}}a)\|_{L^{1}_{T}(L^{p})}+\|\widetilde{R}_{q}\|_{L^{1}_{T}(L^{p})}+\|{\cal P}\D_{q}E_{m}\|_{L^{1}_{T}(L^{p})}\\
&\hspace{1,2cm}+(\frac{\underline{\nu}(p-1)}{p^{2}})^{\frac{1}{m}-1}2^{q(\frac{2}{m}-2)}\big(\|\D_{q}(T_{\n a}\Pi_{2})\|_{L^{1}_{T}(L^{p})}+
2^{q}\|\D_{q}(T_{a}\Pi_{2})\|_{L^{1}_{T}(L^{p})}\\
&\hspace{6,5cm}+\|\D_{q}(T^{'}_{\n \Pi_{2}}a)\|_{L^{1}_{T}(L^{p})}
+
\|{\cal P}g_{q}\|_{L^{m}_{T}(L^{p})}\big),\\
\end{aligned}
\label{11}
\end{equation}
We are now interested by treating the commutator term $\widetilde{R}^{1}_{q}$, we have then by using lemma \ref{alemme3} in the appendix the following estimates with $\alpha<1$:
\begin{equation}
\|\widetilde{R}^{1}_{q}\|_{L^{p}}\lesssim c_{q}\bar{\nu}2^{(-1+\alpha)qs}\|S_{m}a\|_{B^{\frac{N}{p_{1}}+\alpha}
_{p_{1},\infty}}\|Du\|_{B^{s}_{p,r}},
\label{13}
\end{equation}
where $(c_{q})_{q\in\mathbb{Z}}$ is a positive sequence such that $c_{q}\in l^{r}$, and $\bar{\nu}=\mu$.
Note that,
owing to Bernstein inequality, we have:
$$
\begin{aligned}
\|S_{m}a\|_{B^{\frac{N}{p_{1}}+\alpha}_{p_{1},\infty}}&\lesssim
2^{m\alpha}\|a\|_{B^{\frac{N}{p_{1}}}_{p_{1},\infty}}
\end{aligned}
$$
Next we have by corollary \ref{produit2}:
\begin{equation}
\|\widetilde{R}^{2}_{q}\|_{L^{p}}\lesssim c_{q}\bar{\nu}2^{-qs}\|S_{m}a\|_{B^{\frac{N}{p_{1}}}
_{p_{1},r}\cap L^{\infty}}\|u\|_{B^{s+2}_{p,r}},
\label{133}
\end{equation}
Hence, plugging (\ref{13}), (\ref{133}) and (\ref{9}) in (\ref{11}), then multiplying by $2^{qs}$
and summing up on $q\in\mathbb{Z}$ in $l^{r}$, we discover that, for all $t\in[0,T]$:
\begin{equation}
\begin{aligned}
&\|u\|_{\widetilde{L}^{\infty}_{T}(B^{s}_{p,r})}+(\frac{\underline{\nu}(p-1)}{p})^{\frac{1}{m}}
\|u\|_{\widetilde{L}^{m}_{T}(B^{s+\frac{2}{m}}_{p,r})}\leq
\|u_{0}\|_{B^{s}_{p,r}}+
\|{\cal P}f\|_{\widetilde{L}^{1}_{T}(B^{s}_{p,r})}\\
&\hspace{3,5cm}+\|T_{a}\Pi_{1}\|_{\widetilde{L}^{1}_{T}(B^{s+1}_{p,r})}+\|T_{\n a}\Pi_{1}\|_{\widetilde{L}^{1}_{T}(B^{s}_{p,r})}+\|T^{'}_{\n \Pi_{1}}a\|_{\widetilde{L}^{1}_{T}(B^{s}_{p,r})}\\
&+C\bar{\nu}\|a-S_{m}a\|_{\widetilde{L}^{\infty}(
B^{\frac{N}{p_{1}}}_{p_{1},\infty}\cap L^{\infty})}\|u\|_{\widetilde{L}^{1}(B^{s+2}_{p,r})}+2^{m\alpha}\int^{T}_{0}\|a\|_{B^{\frac{N}{p_{1}}}_{p_{1},\infty}}(\tau)\|u\|_{B^{s+1}_{p,r}}(\tau)d\tau\\
&+\big(\frac{\underline{\nu}(p-1)}{p^{2}})^{\frac{1}{m}}
(\|{\cal P}g\|_{\widetilde{L}^{m}_{T}(B^{s-2+\frac{2}{m}}_{p,r})}
+\|T_{a}\Pi_{2}\|_{\widetilde{L}^{m}_{T}(B^{s-1+\frac{2}{m}}_{p,r})}+\|T_{\n a}\Pi_{2}\|_{\widetilde{L}^{m}_{T}(B^{s-2+\frac{2}{m}}_{p,r})}
\\
&\hspace{9cm}+\|T^{'}_{\n \Pi_{2}}a\|_{\widetilde{L}^{m}_{T}(B^{s-2+\frac{2}{m}}_{p,r})}\big),
\end{aligned}
\label{fin1}
\end{equation}
for a constant $C$ depending only on $N$ and $s$.
With our assumption on $\alpha$, $\alpha^{'}$ and $s$, the terms  $\|T_{\n a}\Pi_{1}\|_{\widetilde{L}^{1}_{T}(B^{s}_{p,r})}$ and
$\|T^{'}_{\n \Pi_{1}}a\|_{\widetilde{L}^{1}_{T}(B^{s}_{p,r})}$ may be bounded by:
$$\|a\|_{\widetilde{L}^{\infty}(B^{\frac{N}{p_{1}}+\alpha}_{p,\infty}\cap L^{\infty})}\|\n\Pi_{1}\|
_{\widetilde{L}^{1}_{T}(B^{s-\alpha^{'}}_{p,r})}$$
whereas  $\|T_{\n a}\Pi_{2}\|_{\widetilde{L}^{m}_{T}(B^{s-2+\frac{2}{m}}_{p,r})}$ and $\|T^{'}_{\n \Pi_{2}}a\|_{\widetilde{L}^{m}_{T}(B^{s-2+\frac{2}{m}}_{p,r})}$ may be bounded by:
$$\|a\|_{\widetilde{L}^{\infty}(B^{\frac{N}{p_{1}}+\alpha}_{p,\infty}\cap L^{\infty})}\|\n\Pi_{2}\|
_{\widetilde{L}^{m}_{T}(B^{s-2+\frac{2}{m}-\alpha^{'}}_{p,r})}.$$
Moreover we control $\|T_{a}\Pi_{1}\|_{\widetilde{L}^{1}_{T}(B^{s+1}_{p,r})}$ and $\|T_{a}\Pi_{2}\|_{\widetilde{L}^{m}_{T}(B^{s-1+\frac{2}{m}}_{p,r})}$ by respectively:
$$
\begin{aligned}
&\|a\|_{\widetilde{L}^{\infty}(B^{\frac{N}{p_{1}}}_{p,\infty}\cap L^{\infty})}\|\n\Pi_{1}\|
_{\widetilde{L}^{1}_{T}(B^{s}_{p,r})}.\\
&\|a\|_{\widetilde{L}^{\infty}(B^{\frac{N}{p_{1}}}_{p,\infty}\cap L^{\infty})}\|\n\Pi_{2}\|
_{\widetilde{L}^{m}_{T}(B^{s-2+\frac{2}{m}}_{p,r})}.
\end{aligned}
$$
Hence in view of proposition \ref{propositionA5} and provided that $0<\alpha^{'}<\min(1,\alpha,\frac{s}{2})$ and $s<\frac{N}{p_{1}}$ (which is assumed
in the statement of proposition \ref{linearise}) and $\alpha^{''}\in[0,\alpha^{'}]$,
\begin{equation}
\underline{b}\|\n\Pi_{1}\|_{\widetilde{L}^{1}_{T}(B^{s-\alpha^{''}}_{p,r})}\lesssim{\cal A}_{T}^{\frac{s-\alpha^{''}}{\alpha^{'}}}
\|{\cal Q}f\|_{\widetilde{L}^{1}(B^{s}_{p,r})}.
\label{241}
\end{equation}
On the other hand, by virtue of proposition \ref{produit1}, and of assumption on $\alpha$, $\alpha^{'}$,
$\alpha^{''}$ and $s$, we have:
$$\|{\cal Q}H\|_{\widetilde{L}^{m}_{T}(B^{s-2+\frac{2}{m}-\alpha^{''}}_{p,r})}\lesssim
\|{\cal Q}g\|_{\widetilde{L}^{m}_{T}(B^{s-2+\frac{2}{m}-\alpha^{''}}_{p,r})}+\mu\|a\|_{\widetilde{L}^{\infty}
(B^{\frac{N}{p_{1}}+\alpha}_{p_{1},\infty}\cap L^{\infty})}\|\D u\|_
{\widetilde{L}^{m}_{T}(B^{s-2+\frac{2}{m}-\alpha^{''}}_{p,r})}.$$
As $\alpha^{'}\leq\min(1,\alpha,\frac{1}{2}(s-2+\frac{2}{m}))$, proposition \ref{propositionA5} with $\sigma=s-2+\frac{2}{m}-\alpha^{''}$
(here comes $s>2-\frac{2}{m}$) applies, from which we get for all $\e>0$ ($\e=0$ does if $m\geq2$),
\begin{equation}
\begin{aligned}
&\underline{b}\|\n\Pi_{2}\|_{\widetilde{L}^{m}_{T}(B^{s-2+\frac{2}{m}-\alpha^{''}}_{p,r})}
\lesssim{\cal A}_{T}^{\frac{s-2+\frac{2}{m}+\e}{\alpha^{'}}}\big(
\|{\cal Q}g\|_{\widetilde{L}^{m}_{T}(B^{s-2+\frac{2}{m}-\alpha^{''}}_{p,r}}\\
&\hspace{5cm}+\mu\|a\|_{\widetilde{L}^{\infty}
(B^{\frac{N}{p_{1}}+\alpha}_{p_{1},\infty}\cap L^{\infty})}\|\D u\|_{\widetilde{L}^{m}_{T}(B^{s-2+\frac{2}{m}
-\alpha^{''}}_{p,r})}\big).
\end{aligned}
\label{251}
\end{equation}
Let $X(t)=\|u\|_{L^{\infty}_{t}(B^{s}_{p,r})}+\nu \underline{b}\|u\|_{L^{1}_{t}(B^{s+2}_{p,r})}$.
Assuming that $m$ has been chosen so
large as to satisfy:
$$C\bar{\nu}\|a-S_{m}a\|_{L^{\infty}_{T}(B^{\frac{N}{p_{1}}}_{p_{1},\infty}\cap L^{\infty})}\leq\underline{\nu},$$
and  by interpolation we get:
\begin{equation}
C\bar{\nu}2^{m\alpha}\|a\|_{B^{\frac{N}{p_{1}}}_{p_{1},\infty}}\|u\|_{B^{s+2}_{p,r}}\leq\kappa\underline{\nu}+\frac{C^{2}
\bar{\nu}^{2}2^{2m\alpha}}
{4\kappa\underline{\nu}}
\|a\|^{2}_{B^{\frac{N}{p_{1}}}_{p_{1},\infty}}\|u\|_{B^{s}_{p,r}},
\label{2511}
\end{equation}
Plugging (\ref{241}), (\ref{251}) and (\ref{2511}) in (\ref{fin1}), we end up with:
$$
\begin{aligned}
&X(T)\leq\|u_{0}\|_{B^{s}_{p,r}}+{\cal A}_{T}^{\frac{s}{\alpha^{'}}}\biggl(
\|{\cal P}f\|_{\widetilde{L}^{1}_{t}(B^{s}_{p,r})}+\|{\cal P}g\|_{\widetilde{L}^{m}_{t}(B^{s-2+\frac{2}{m}}_{p,r})}+C\int^{t}_{0}\big(\frac{\bar{\nu}^{2}}
{\underline{\nu}}2^{2m\alpha}
\|a\|^{2}_{B^{\frac{N}{p_{1}}}_{p_{1},\infty}}(\tau)\\
&\hspace{7cm}\times X(\tau)\big)d\tau+(\frac{\underline{\nu}(p-1)}{p})^{\frac{1}{m}}{\cal A}_{T}\|u\|_{\widetilde{L}^{1}_{T}(B^{s+2-\alpha^{'}}_{p,r})}\biggl).
\end{aligned}
$$
Gr\"onwall lemma then leads to the desired inequality.
\subsubsection*{Case $p=2$}
In this case we don't need of condition of smallness on $\|a\|_{B^{\frac{N}{p_{1}}}_{p_{1},\infty}\cap L^{\infty}}$, indeed
the bad terms as $\widetilde{R}^{2}_{q}$ or $2^{q}\|\D_{q}(T_{a}\Pi)\|_{L^{2}}$ disapear in the integration by parts as
${\rm div}u_{q}=0$. So we can follow the same procedure and conclude.
{\hfill $\Box$}
\section{The mass conservation equation}
\label{section4}
\subsection{Losing estimates for transport equation}
We now focus on the mass equation associated to $(\ref{NSC})$:
\begin{equation}
\begin{cases}
\begin{aligned}
&\p_{t}a+v\cdot\n a=g,\\
&a_{/t=0}=a_{0}.
\end{aligned}
\end{cases}
\label{21}
\end{equation}
We will precise in the sequel the regularity of $a_{0}$, $v$ and $g$. In this section we intend to recall some result on transport equation associated to vector fields which are not
Lipschitz with respect to the space variable. Since we still have in mind to get regularity theorems, those vector
field cannot be to rough.
In order to measure precisely the regularity of the vector field $v$, we shall introduce the following notation:
\begin{equation}
V^{'}_{p_{1},\alpha}(t)=\sup_{j\geq 0}\frac{2^{j\frac{N}{p_{1}}}\|\n S_{j}v(t)\|_{L^{p_{1}}}}{(j+1)^{\alpha}}<+\infty.
\label{331}
\end{equation}
Let us remark that if $p_{1}=+\infty$ then $V^{'}_{p_{1},\alpha}$ is exactly the norm $\|B_{\Gamma}\|$ of definition \ref{defChemin1}.
\subsubsection{Limited loss of regularity}
In this section, we make the assumption that there exists some $\alpha\in]0,1[$ such that the function $V^{'}_{p_{1},\alpha}$ defined in (\ref{331}) be locally integrable. We will show that in the case $\alpha=1$, then a linear loss of regularity may occur.
In the theorem below, Bahouri, Chemin and Danchin show in \cite{BCD} that if $\alpha\in]0,1[$
then the loss of regularity in the estimate is arbitrarily small.
\begin{theorem}
\label{transport1}
Let $(p,p_{1})$ be in $[1,+\infty]^{2}$ such that $1\leq p\leq p_{1}$ and $\sigma$ satisfying
$\sigma>-1-N\min (\frac{1}{p_{1}},\frac{1}{p^{'}})$. Assume that $\sigma<1+\frac{N}{p_{1}}$ and that
$V^{'}_{p_{1},\alpha}\in]0,1[$ is in $L^{1}([0,T])$. Let $a_{0}\in B^{\sigma}_{p,\infty}$ and $g\in \widetilde{L}_{T}^{1}(B^{\sigma}_{p,\infty})$. Then the equation
(\ref{21}) has a unique solution $a\in C([0,T],\cap_{\sigma^{'}<\sigma}B^{\sigma^{'}}_{p,\infty})$ and the
following estimate
holds for all small enough $\e$:
$$\|a\|_{\widetilde{L}^{\infty}_{T}(B^{\sigma-\e}_{p,\infty})}\leq C\big(\|a_{0}\|_{B^{\sigma}_{p,\infty}}+\|g\|_{\widetilde{L}_{T}^{1}(B^{\sigma}_{p,\infty})}\big)\exp\big(
\frac{C}{\e^{\frac{\alpha}{1-\alpha}}}(V_{p_{1},\alpha}(T))^{\frac{1}{1-\alpha}}\big),$$
where $C$ depends only on $\alpha$, $p$, $p_{1}$, $\sigma$ and $N$.
\end{theorem}
In the following proposition, we are interested in showing a control of the high frequencies on the density when
$u$ is not Lipschitz. Indeed we recall that in the proposition \ref{linearise} when $p=2$, we need to control the high frequencies of the density. In particular the following proposition is useful only in the case of theorem \ref{theo11}.
\begin{proposition}
\label{proptransport1}
Let $(p,p_{1})$ be in $[1,+\infty]^{2}$ such that $1\leq p\leq p_{1}$ and $\sigma$ satisfying
$\sigma>-1-N\min (\frac{1}{p_{1}},\frac{1}{p^{'}})$. Assume that $\sigma<1+\frac{N}{p_{1}}$ and that
$V^{'}_{p_{1},\alpha}\in]0,1[$ is in $L^{1}([0,T])$. Let $a_{0}\in B^{\sigma}_{p,\infty}$ and $g\in \widetilde{L}_{T}^{1}(B^{\sigma}_{p,\infty})$, the equation
(\ref{21}) has a unique solution $a\in C([0,T],\cap_{\sigma^{'}<\sigma}B^{\sigma^{'}}_{p,\infty})$ and the
following estimate
holds for all small enough $\e$:
$$
\begin{aligned}
&\sum_{l\geq m}2^{(\sigma-\e)l}\|\D_{l}a(t^{'})\|_{L^{\infty}(L^{p})}
\lesssim \sum_{l\geq m}( 2^{\sigma l}\|\D_{l}a_{0}\|_{L^{p}})+C\eta^{\frac{\alpha}{1-\alpha}}
\int^{t}_{0}V^{'}_{p_{1},\alpha}(t^{'})\\
&\hspace{4cm}\times(\|a_{0}\|_{B^{\sigma}_{p,\infty}}+\|g\|_{\widetilde{L}_{T}^{1}(B^{\sigma}_{p,\infty})})
\exp\big(\frac{C}{\e^{\frac{\alpha}{1-\alpha}}}
(V_{p_{1},\alpha}(t^{'}))^{\frac{1}{1-\alpha}}\big)dt^{'}.
\end{aligned}
$$
where $C$ depends only on $\alpha$, $p$, $p_{1}$, $\sigma$ and $N$.
\end{proposition}
{\bf Proof:}\\
\\
By using the proof of Bahouri, Chemin and Danchin in \cite{BCD} one can write:
\begin{equation}
2^{(2+j)\sigma_{t}}\|\D_{l}a(t)\|_{L^{p}}\leq 2^{(2+l)\sigma}\|\D_{l}f_{0}\|_{L^{p}}+
C(\frac{2C}{\eta\log2})^{\frac{\alpha}{1-\alpha}}\int^{t}_{0}V^{'}_{p_{1},\alpha}(t^{'})\|a(t^{'})\|_{B^
{\sigma_{t^{'}}}_{p,\infty}}dt^{'}.
\label{3.36}
\end{equation}
Whence taking the supremum over $l\geq m$, we get
$$ \sup_{t^{'}\in[0,t]}\,\sup_{l\geq m}( 2^{\sigma_{t^{'}}l}\|\D_{l}a(t^{'})\|_{L^{p}})
\lesssim \sup_{l\geq m}( 2^{\sigma l}\|\D_{l}a_{0}\|_{L^{p}})+C\eta^{\frac{\alpha}{1-\alpha}}
\int^{t}_{0}V^{'}_{p_{1},\alpha}(t^{'})\|a(t^{'})\|_{B^{\sigma_{t^{'}}}_{p,\infty}}dt^{'}.
$$ We apply Gr\"onwall inequality in (\ref{3.36}) and we insert in previous inequality:
$$
\begin{aligned}
&\sup_{t^{'}\in[0,t]}\,\sup_{l\geq m}( 2^{\sigma_{t^{'}}l}\|\D_{l}a(t^{'})\|_{L^{p}})
\lesssim \sup_{l\geq m}( 2^{\sigma l}\|\D_{l}a_{0}\|_{L^{p}})+C\eta^{\frac{\alpha}{1-\alpha}}
\|f_{0}\|_{B^{\sigma}_{p,\infty}}\int^{t}_{0}V^{'}_{p_{1},\alpha}(t^{'})\\
&\hspace{8,5cm}\times
\exp\big(\frac{C}{\e^{\frac{\alpha}{1-\alpha}}}
(V_{p_{1},\alpha}(t^{'}))^{\frac{1}{1-\alpha}}\big)dt^{'}.
\end{aligned}
$$
which leads to the proposition.
\hfill {$\Box$}
\begin{remarka}
In the sequel, we will use the theorem \ref{transport1} and the proposition \ref{proptransport1} when $p_{1}=\infty$ and $\alpha=1+\e^{'}-\frac{1}{r}$.
Indeed we will have $u$ is in $\widetilde{L}^{1}(B^{\frac{N}{p_{2}}+1}_{p_{2},r})$ (with $1\leq r<\infty$) and according proposition
\ref{propChemin2} and definition \ref{defChemin1} we get:
$$\int^{t}_{0}V^{'}_{\infty,\alpha}(t^{'})dt^{'}\lesssim\|u\|_{\widetilde{L}^{1}(B^{\frac{N}{p_{2}}+1}_{p_{2},r})}.$$
So it will allow to get estimates on the density with an arbitrarly small loss of regularity.
\end{remarka}
\subsubsection{Linear loss of regularity in Besov spaces}
Tis section is devoted to the estimates with linear loss of regularity. Remind that $u$ is log-lipschitz only if only if there some constant $C$ such that:
$$\|\n S_{j}\|_{L^{\infty}}\leq C(j+1)\;\;\;\mbox{for all}\:\:j\geq-1.$$
So Bahouri, Chemin, Danchin show in \cite{BCD} the following theorem.
\begin{theorem}
\label{transport2}
Let $1\leq p\leq p_{1}\leq\infty$ and $s_{1}\in\R$ satisfies $\sigma>-N\min (\frac{1}{p_{1}},\frac{1}{p^{'}})$.
Let $\sigma$ in $]s_{1},1+\frac{N}{p_{1}}[$ and v a vector field. There exists a constant $C$ depending only on
$p$, $p_{1}$, $\sigma$, $s_{1}$ and $N$ such that for any $\lambda>C,\,T>0$ and any nonnegative integrable function
$W$ over $[0,T]$ such that if $\sigma_{T}\geq s_{1}$ with:
$$\sigma_{t}=\sigma-\lambda\int^{t}_{0}\big(V^{'}_{p_{1},1}(t^{'})+W(t^{'})\big)dt^{'}$$
then the following property holds true.\\
Let $a_{0}\in B^{\sigma}_{p,\infty}$ and $g=g_{1}+g_{2}$ with, for all $t\in[0,T]$, $g_{1}(t)\in B^{\sigma_
{t}}_{p,\infty}$ and:
$$\forall j\geq-1,\;\|\D_{j}g_{2}(t)\|_{L^{p}}\leq 2^{-j\sigma_{t}}(2+j)W(t)\|f(t)\|_{B^{\sigma_{t}}_{p,\infty}}.$$
Let $a\in {\cal C}([0,T];B^{s_{1}}_{p,\infty})$ be a solution of (\ref{21}). Then the following estimate holds:
$$\sup_{t\in[0,T]}\|a(t)\|_{B^{\sigma_{t}}_{p,\infty}}\leq\frac{\lambda}{\lambda-C}\big(\|a_{0}\|_{B^{\sigma}
_{p,\infty}}+
\int^{T}_{0}\|g_{1}(t)\|_{B^{\sigma}_{p,\infty}}dt\big).$$
\end{theorem}
In the following proposition, we generalize this result to the high frequencies of the density $\rho$.
\begin{proposition}
\label{proptransport2}
Let $1\leq p\leq p_{1}\leq\infty$ and $s_{1}\in\R$ satisfies $\sigma>-N\min (\frac{1}{p_{1}},\frac{1}{p^{'}})$.
Let $\sigma$ in $]s_{1},1+\frac{N}{p_{1}}[$ and v a vector field. There exists a constant $C$ depending only on
$p$, $p_{1}$, $\sigma$, $s_{1}$ and $N$ such that for any $\lambda>C,\,T>0$ and any nonnegative integrable function
$W$ over $[0,T]$ such that if $\sigma_{T}\geq s_{1}$ with:
$$\sigma_{t}=\sigma-\lambda\int^{t}_{0}\big(V^{'}_{p_{1},1}(t^{'})+W(t^{'})\big)dt^{'}$$
then the following property holds true.\\
Let $a_{0}\in B^{\sigma}_{p,\infty}$ and $g=g_{1}+g_{2}$ with, for all $t\in[0,T]$, $g_{1}(t)\in B^{\sigma_
{t}}_{p,\infty}$ and:
$$\forall j\geq-1,\;\|\D_{j}g_{2}(t)\|_{L^{p}}\leq 2^{-j\sigma_{t}}(2+j)W(t)\|f(t)\|_{B^{\sigma_{t}}_{p,\infty}}.$$
Let $a\in {\cal C}([0,T];B^{s_{1}}_{p,\infty})$ be a solution of (\ref{21}). Then the following estimate holds:
$$\sup_{t\in{0,T}}\sum_{l\geq m}2^{l\sigma_{t}}\|\D_{l}a\|_{L^{p}}\leq\frac{\lambda}{\lambda-C}\big(\|a_{0}\|_{B^{\sigma}
_{p,\infty}}+
\int^{T}_{0}\|g_{1}(t)\|_{B^{\sigma}_{p,\infty}}dt\big).$$
\end{proposition}
{\bf Proof:}
\\
\\
The proof follows the same line as in proposition \ref{proptransport1}.
\hfill {$\Box$}
\begin{remarka}
In the sequel, we will use the theorem \ref{transport2} and the proposition \ref{proptransport2}
when $p_{1}=\infty$.
Indeed we will have $u$ in $\widetilde{L}^{1}(B^{\frac{N}{p_{2}}+1}_{p_{2},\infty})$,
then according proposition \ref{derivation,interpolation}, $\n u$ is in
$\widetilde{L}^{1}(B^{0}_{\infty,\infty})$so we obtain the following control:
$$\int^{t}_{0}\sup_{l\geq0}\big(\frac{\|\n S_{l}u(t^{'},\cdot)\|_{L^{\infty}}}{l+1}\big)dt^{'}\lesssim\|u\|_{\widetilde{L}^{1}
(B^{\frac{N}{p_{2}}+1}_{p_{2},\infty})}.$$
So it will allow to get estimates on the density with a small loss of regularity on a small time interval.
\end{remarka}
\section{Proof of the existence for theorem \ref{theo1}}
\label{section5}
We use a standard scheme:
\begin{enumerate}
\item We smooth out the data and get a sequence of smooth solutions
$(a^{n},u^{n})_{n\in\mathbb{N}}$ to (\ref{0.6}) on a bounded interval $[0,T^{n}]$ which may depend on $n$.
\item We exhibit a positive lower bound $T$ for $T^{n}$, and prove uniform estimates in the space:
$$E_{T}=\widetilde{C}_{T}(B^{\frac{N}{p_{1}}+\frac{\e}{2}}_{p_{1},\infty})\times\big(
\widetilde{C}_{T}(B^{\frac{N}{p_{2}}-1}_{p_{2},r})\cap\widetilde{L}^{1}(B^{\frac{N}{p_{�}}+1}_{p_{2}
,r}))\big)^{N}.$$
for the smooth solution $(a^{n},u^{n})$.
\item We use compactness to prove that the sequence converges, up to extraction, to a solution of (\ref{0.6}).
\end{enumerate}
\subsubsection*{Construction of approximate solutions}
We smooth on the data as follows:
$$a_{0}^{n}=S_{n}a_{0},\;\;u_{0}^{n}=S_{n}u_{0}\;\;\;\mbox{and}\;\;\;f^{n}=S_{n}f.$$
Note that we have:
$$\forall l\in\mathbb{Z},\;\;\|\D_{l}a^{n}_{0}\|_{L^{p_{1}}}\leq\|\D_{l}a_{0}\|_{L^{p_{1}}}\;\;\;\mbox{and}\;\;\;\|a^{n}_{0}\|
_{B^{\frac{N}{p_{1}}+\e^{'}}_{p,\infty}}\leq \|a_{0}\|_{B^{\frac{N}{p_{1}}+\e^{'}}_{p,\infty}},$$
and similar properties for $u_{0}^{n}$ and $f^{n}$, a fact which will be used repeatedly during the next
steps. Now, according \cite{AP}, one can solve (\ref{0.6}) with the smooth data $(a_{0}^{n},u_{0}^{n},f^{n})$.
We get a solution $(a^{n},u^{n})$ on a non trivial time interval $[0,T_{n}]$ such that:
\begin{equation}
\begin{aligned}
&a^{n}\in\widetilde{C}([0,T_{n}),H^{N+\e}),\;\;u^{n}\in{\cal C}([0,T_{n}),H^{N})\cap\widetilde{L}^{1}_{T_{n}}
(H^{N+2})\\
&\hspace{8cm}\;\;\mbox{and}\;\;\n\Pi^{n}\in\widetilde{L}^{1}_{T_{n}}(H^{N}).
\end{aligned}
\label{a26}
\end{equation}
\subsubsection*{Uniform bounds}
Let $T_{n}$ be the lifespan of $(a_{n},u_{n})$, that is the supremum of all $T>0$ such that (\ref{0.6}) with initial data
$(a_{0}^{n},u_{0}^{n})$ has a solution which satisfies (\ref{a26}). Let $T$ be in $(0,T_{n})$,
we aim at getting uniform estimates in $E_{T}$ for $T$ small enough. For that, we need to introduce the solution
$(u^{n}_{L},\n\Pi_{L}^{n})$ to the nonstationary Stokes system:
$$
(L)\;\;\;\;
\begin{cases}
\begin{aligned}
&\p_{t}u_{L}^{n}-\mu\D u_{L}^{n}+\n\Pi^{n}_{L}=f^{n},\\
&{\rm div}u_{L}^{n}=0,\\
&(u^{n}_{L})_{t=0} =u^{n}_{0}.
\end{aligned}
\end{cases} $$
Now, the vectorfields $\tilde{u}^{n}=u^{n}-u^{n}_{L}$ and $\n\Pi^{n}=\n\Pi^{n}_{L}+\n\tilde{\Pi}^{n}$
satisfy the parabolic system:
\begin{equation}
\begin{cases}
\begin{aligned}
&\p_{t}\tilde{u}^{n}-\mu(1+a^{n})\D\tilde{u}^{n}+(1+a^{n})\n\tilde{\Pi}^{n}=H(a^{n},u^{n},\n\Pi^{n}),\\
&{\rm div}\tilde{u}^{n}=0,\\
&\tilde{u}^{n}(0)=0,
\end{aligned}
\end{cases}
\label{3momentum}
\end{equation}
which has been studied in proposition \ref{linearise} with:
$$H(a^{n},u^{n},\n\Pi^{n})=a^{n}(\mu\D u^{n}_{L}-\n\Pi_{L}^{n})-u^{n}\cdot\n u^{n}.$$
Define $m\in\mathbb{Z}$ by:
\begin{equation}
m=\inf\{ p\in\mathbb{Z}/\;\; 2\bar{\nu}\sum_{\l\geq q}
2^{l\frac{N}{p_{1}}}\|\D_{l}a_{0}\|_{L^{p_{1}}}\leq c\bar{\nu}\}
\label{def}
\end{equation}
where $c$ is small enough positive constant (depending only $N$) to be fixed hereafter.\\
Let:
$$\bar{b}=1+\sup_{x\in\R^{N}}a_{0}(x),\;A_{0}=1+2\|a_{0}\|_{B^{\frac{N}{p_{1}}+\e^{'}}_{p_{1},\infty}},\;
U_{0}=\|u_{0}\|_{B^{\frac{N}{p_{2}}-1}_{p_{2},r}}+
\|f\|_{\widetilde{L}^{1}(B^{\frac{N}{p_{2}}-1}_{p_{2},r})},$$
and $\widetilde{U}_{0}=2CU_{0}+4C\bar{\nu}A_{0}$ (where $C$ stands for a large enough constant depending only $N$ which will be determined when applying proposition \ref{produit1}, \ref{linearise} and \ref{proptransport1} in the following computations.) We assume that the following inequalities are fulfilled for some $\eta>0$, $\alpha>0$:
$$
\begin{aligned}
&({\cal H}_{1})\hspace{2cm}&\|a^{n}-S_{m}a^{n}\|_{\widetilde{L}^{\infty}_{T}(B^{\frac{N}{p_{1}}}_{p_{1},\infty})\cap L^{\infty}}\leq c\underline{\nu}\bar{\nu}^{-1},\\
&({\cal H}_{2})\hspace{1cm}&C\bar{\nu}^{2}T\|a^{n}\|^{2}_{\widetilde{L}^{\infty}_{T}(B^{\frac{N}{p_{1}}}_{p_{1},\infty})\cap L^{\infty}}\leq 2^{-2m}\underline{\nu},\\
&({\cal H}_{3})\hspace{2cm}&\frac{1}{2}\underline{b}\leq 1+a^{n}(t,x)\leq 2\bar{b}\;\;\mbox{for all}\;\;(t,x)\in[0,T]\times\R^{N},\\
&({\cal H}_{4})\hspace{2cm}&\|a^{n}\|_{\widetilde{L}^{\infty}_{T}(B^{\frac{N}{p_{1}}+\frac{\e}{2}}_{p_{1},\infty})\cap L^{\infty}}\leq A_{0},\\
&({\cal H}_{5})\hspace{2cm}&\|u^{n}_{L}\|_{\widetilde{L}^{1}(B^{\frac{N}{p_{2}}+1}_{p_{2},r})}\leq \eta,\\
&({\cal H}_{6})\hspace{1,5cm}&\|\widetilde{u}^{n}\|_{\widetilde{L}^{\infty}(B^{\frac{N}{p_{2}}-1}_{p_{2},r})}
+\underline{\nu}
\|\widetilde{u}^{n}\|_{\widetilde{L}^{1}(B^{\frac{N}{p_{2}}+1}_{p_{2},r})}\leq \widetilde{U}_{0}\eta,\\
&({\cal H}_{7})\hspace{1,5cm}&\|\n\Pi^{n}_{L}\|_{\widetilde{L}^{1}(B^{\frac{N}{p_{2}}-1}_{p_{2},r})}
\leq \eta,\\
&({\cal H}_{8})\hspace{1,5cm}&\|\n\widetilde{\Pi}^{n}\|_{\widetilde{L}^{1}(B^{\frac{N}{p_{2}}-1}_{p_{2},r})}
\leq \widetilde{\Pi}_{0}\eta,
\end{aligned}
$$
Remark that since:
$$1+S_{m}a^{n}=1+a^{n}+(S_{m}a^{n}-a^{n}),$$
assumptions $({\cal H}_{1})$ and $({\cal H}_{3})$ insure that:
\begin{equation}
\inf_{(t,x)\in[0,T]\times\R^{N}}(1+S_{m}a^{n})(t,x)\geq\frac{1}{4}\underline{b},
\label{inemin}
\end{equation}
provided $c$ has been chosen small enough (note that $\frac{\underline{\nu}}{\bar{\nu}}\leq\bar{b}$).\\
We are going to prove that under suitable assumptions on $T$ and $\eta$ (to be specified below) if condition $({\cal H}_{1})$ to $({\cal H}_{8})$ are satisfied, then they are actually satisfied with strict inequalities. Since all those conditions depend continuously on the time variable and
are strictly satisfied initially, a basic boobstrap argument insures that $({\cal H}_{1})$ to $({\cal H}_{8})$ are indeed satisfied for a small $T^{'}_{n}\leq T^{n}$.
First we shall assume that $\eta$ satisfies:
\begin{equation}
C(1+\underline{\nu}^{-1}\widetilde{U}_{0})\eta\leq\log 2
\label{1conduti}
\end{equation}
so that denoting:
$$(\tilde{V}^{n})_{p_{2},1-\frac{1}{r}}(t)=\int^{t}_{0}(\tilde{V}^{n})_{p_{2},1-\frac{1}{r}}^{'}(s)ds\;\;\;\mbox{and}\;\;
(V^{n}_{L})_{p_{2},1-\frac{1}{r}}(t)=\int^{t}_{0}(V_{L}^{n})_{p_{2},1-\frac{1}{r}}^{'}(s)ds,$$
with:
$$
\begin{aligned}
&(\tilde{V}^{n})_{p_{2},1-\frac{1}{r}}^{'}(s)=\sup_{l\geq0}\big(\frac{2^{l\frac{N}{p_{2}}}\|\n S_{l}
\tilde{u}^{n}(s)\|_{L^{p_{2}}}}{(l+1)^{1-\frac{1}{r}}}\big)ds\;\;\;\mbox{and}\\
&\hspace{4cm}(V^{n}_{L})_{p_{2},1-\frac{1}{r}}^{'}(s)=\sup_{l\geq0}\big(\frac{2^{l\frac{N}{p_{2}}}\|\n S_{l}
u_{L}^{n}(s)\|_{L^{p_{2}}}}{(l+1)^{1-\frac{1}{r}}}\big)ds,
\end{aligned}
$$
We recall now that according proposition \ref{propChemin2}:
$$
\begin{aligned}
&(\tilde{V}^{n})_{p_{2},1-\frac{1}{r}}(t)\leq C\|\tilde{u}^{n}\|_{\widetilde{L}^{1}_{t}
(B^{\frac{N}{p_{2}}+1}_{p_{2},r})}\;\;\;
\mbox{and}\;\;\;(V_{L}^{n})_{p_{2},1-\frac{1}{r}}(t)\leq C\|u^{n}_{L}\|_{\widetilde{L}^{1}_{t}(B^{\frac{N}{p_{2}}+1}_{p_{2},r})},
\end{aligned}
$$
we have, according to $({\cal H}_{5})$ and $({\cal H}_{6})$:
\begin{equation}
\begin{aligned}
&e^{C (\frac{e}{2})^{1-r}\big(\|\tilde{u}^{n}\|_{\widetilde{L}^{1}_{t}(B^{\frac{N}{p_{2}}+1}
_{p_{2},r})}+
\|u^{n}_{L}\|_{\widetilde{L}^{1}_{t}(B^{\frac{N}{p_{2}}+1}_{p_{2},r})}\big)}<2.
\end{aligned}
\label{1ineimpca}
\end{equation}
In order to bound $a^{n}$ in $\widetilde{L}^{\infty}_{T}(B^{\frac{N}{p_{1}}+\frac{\e}{2}}_{p_{1},\infty})$, we apply
theorem \ref{transport1} and get:
\begin{equation}
\begin{aligned}
&\|a^{n}\|_{\widetilde{L}^{\infty}_{T}(B^{\frac{N}{p_{1}}+\frac{\e}{2}}_{p_{1},\infty})}
< e^{C(\frac{\e}{2})^{1-r}}\big(\|\tilde{u}^{n}\|_{\widetilde{L}^{1}_{t}(B^{\frac{N}{p_{2}}}_
{p_{2},r})}+
\|u^{n}_{L}\|_{\widetilde{L}^{1}_{t}(B^{\frac{N}{p_{2}}}_{p_{2},r})}\big)
\|a_{0}\|_{B^{\frac{N}{p_{1}}+\e}_{p_{1},\infty}}.
\end{aligned}
\label{inetranspr}
\end{equation}
Moreover as we know that $\|a^{n}\|_{L^{\infty}}\leq \|a_{0}\|_{L^{\infty}}$,  $({\cal H}_{4})$ is satisfied with a strict inequality.
Next, applying classical proposition on heat equation ( see \cite{Da1}) yields:
\begin{equation}
\|u_{L}^{n}\|_{\widetilde{L}^{\infty}_{T}(B^{\frac{N}{p_{2}}-1}_{p_{2},r})}\leq U_{0},
\label{34}
\end{equation}
\begin{equation}
\kappa\nu\|u^{n}_{L}\|_{\widetilde{L}^{1}_{T}(B^{\frac{N}{p_{2}}+1}_{p_{2},r})}
\leq\big(\sum_{l\in\mathbb{Z}}2^{lr(\frac{N}{p_{2}}-1)}(1-e^{-\kappa\nu2^{2l}T})^{r}
(\|\D_{l}u_{0}\|^{r}_{L^{p_{2}}}+\|\D_{l}f\|^{r}_{L^{1}(\R^{+},L^{p_{2}})})\big)^{\frac{1}{r}}.
\label{35}
\end{equation}
Hence taking $T$ such that:
\begin{equation}
\big(\sum_{l\in\mathbb{Z}}2^{lr(\frac{N}{p_{2}}-1)}(1-e^{-\kappa\nu2^{2l}T})^{r}(\|\D_{l}u_{0}\|^{r}_{L^{2}}
+\|\D_{l}f\|^{r}_{L^{1}(\R^{+},L^{2})})\big)^{\frac{1}{r}}<\kappa\eta\nu,
\label{36}
\end{equation}
insures that $({\cal H}_{5})$ is strictly verified.\\
Since $({\cal H}_{1})$, $({\cal H}_{2})$ and (\ref{inemin}) are satisfied,
proposition \ref{linearise} may be applied with $\alpha=\frac{\e}{2}$. We get:
$$
\begin{aligned}
&\|\widetilde{u}^{n}\|_{\widetilde{L}^{\infty}_{T}(B^{\frac{N}{p_{2}}-1}_{p_{2},r})}+\underline{\nu}
\|\widetilde{u}^{n}\|_{\widetilde{L}^{1}_{T}(B^{\frac{N}{p_{2}}+1}_{p_{2},r})}\leq Ce^{CZ^{n}_{m}(T)}\biggl(
\|u_{0}\|_{B^{\frac{N}{p_{2}}-1}_{p_{2},r}}+{\cal A}_{T,n}^{\kappa}\times\\
&\big(\|a^{n}(\D u^{n}_{L}-\n\Pi_{L}^{n})\|_{\widetilde{L}^{1}(B^{\frac{N}{p_{2}}-1}_{p_{2},r})}
+\|u^{n}\cdot\n u^{n}\|_{\widetilde{L}^{1}(B^{\frac{N}{p_{2}}-1}_{p_{2},r})}+{\cal A}_{T,n}\|u\|_
{\widetilde{L}^{1}_{T}
(B^{\frac{N}{p_{2}}+1-\alpha^{'}}_{p_{2},r})}\big)\biggl).
\end{aligned}
$$
with $Z^{n}_{m}(T)=2^{m\e}\bar{\nu}^{2}\underline{\nu}^{-1}\int^{T}_{0}\|a^{n}\|^{2}_{B^{\frac{N}{p_{1}}}_
{p_{1},\infty}\cap L^{\infty}}d\tau.$
Next using Bony's decomposition and ${\rm div}u^{n}=0$, one can write:
$${\rm div}(u^{n}\otimes u^{n})=T_{\p_{j}v}(u^{n})^{j}+T_{(u^{n})^{j}}\p_{j}u^{n}+\p_{j}R(u^{n},(u^{n})^{j}),$$
with the summation convention over repeated indices.\\
Hence combining proposition 1.4.1 and 1.4.2 in \cite{Da1} with the fact that
$\widetilde{L}^{\rho}_{T}(B^{\frac{N}{p_{2}}-\frac{1}{2}}_{p_{2},r})\hookrightarrow
\widetilde{L}^{\rho}_{T}(B^{-\frac{1}{2}}_{\infty,\infty})$ for $\rho=\frac{4}{3}$ or $\rho=4$, we get:
$$\|{\rm div}(u^{n}\otimes u^{n})\|_{\widetilde{L}^{1}_{T}(B^{\frac{N}{p_{2}}-1}_{p_{2},r})}\leq
C\|u^{n}\|_{\widetilde{L}^{\frac{4}{3}}_{T}(B^{\frac{N}{p_{2}}+\frac{1}{2}}_{p_{2},r})}
\|u^{n}\|_{\widetilde{L}^{4}_{T}(B^{\frac{N}{p_{2}}-\frac{1}{2}}_{p_{2},r})}.$$
By taking advantage of proposition \ref{linearise}, \ref{produit1}, \ref{derivation,interpolation} and Young'inequality, we end up with:
$$
\begin{aligned}
&\|\widetilde{u}^{n}\|_{\widetilde{L}^{\infty}_{T}(B^{\frac{N}{p_{2}}-1}_{p_{2},r})}+\underline{\nu}
\|\widetilde{u}^{n}\|_{\widetilde{L}^{1}_{T}(B^{\frac{N}{p_{2}}+1}_{p_{2},r})}\leq e^{CZ^{n}_{m}(T)}\biggl(
\|u_{0}\|_{B^{\frac{N}{p_{2}}-1}_{p_{2},r}}+
\big(
\|a^{n}\|_{\widetilde{L}^{\infty}_{T}(B^{\frac{N}{p_{1}}+\frac{\e}{2}}_{p_{1},\infty}\cap L^{\infty})}+1\big)^{\kappa}\\
&\times\big(
C\|u^{n}_{L}\|_{\widetilde{L}^{1}_{T}(B^{\frac{N}{p_{2}}+1}_{p_{2},r})}(
\bar{\nu}\|a^{n}\|_{L^{\infty}_{T}(B^{\frac{N}{p_{1}}}_{p_{1},\infty}\cap L^{\infty})}
+\|u^{n}_{L}\|_{\widetilde{L}^{\infty}_{T}(B^{\frac{N}{p_{2}}-1}_{p_{2},r})})+\|\n\Pi^{n}_{L}\|_{\widetilde{L}^{1}_{T}(B^{
\frac{N}{p_{2}}-1}_{p_{2},r})}\times\\
&\|a^{n}\|_{L^{\infty}_{T}(B^{\frac{N}{p_{1}}}_{p_{1},\infty}\cap L^{\infty})}\big)+(\|a^{n}\|_{\widetilde{L}^{\infty}_{T}(B^{\frac{N}{p_{1}}+\frac{\e}{2}}_{p_{1},\infty}\cap L^{\infty})}+1)T^{\frac{\e}{2}}\|\widetilde{u}^{n}\|_{\widetilde{L}^{1}_{T}(B^{\frac{N}{p_{2}}+1}_{p_{2},r})}\biggl).
\end{aligned}
$$
with $C=C(N)$ and $C_{g}=(N,g,\underline{b},\bar{b})$. Now, using assumptions $({\cal H}_{2})$, $({\cal H}_{4})$, $({\cal H}_{5})$, $({\cal H}_{6})$ and
$({\cal H}_{7})$, and inserting (\ref{1ineimpca}) in the previous inequality and choosing $T$ enough small gives:
$$\|\widetilde{u}^{n}\|_{\widetilde{L}^{\infty}_{T}(B^{\frac{N}{p_{2}}-1}_{p_{2},r})}+
\|\widetilde{u}^{n}\|_{\widetilde{L}^{1}_{T}(B^{\frac{N}{p_{2}}+1}_{p_{2},r})}\leq2C(\bar{\nu}A_{0}+U_{0})\eta+2C_{g}TA_{0},$$
hence $({\cal H}_{6})$ is satisfied with a strict inequality provided:
\begin{equation}
C_{g}T<C\bar{\nu}\eta.
\label{38}
\end{equation}
To show that $({\cal H}_{7})$ and  $({\cal H}_{8})$ are strictly verified on $(0,T^{'}_{n})$, we proceed similarly as for
$({\cal H}_{5})$ and  $({\cal H}_{6})$.
We now have to check whether $({\cal H}_{1})$ is satisfied with strict inequality. For that we apply proposition \ref{proptransport1} which yields for all $m\in\mathbb{Z}$,
\begin{equation}
\sup_{l\geq m}2^{l\frac{N}{p_{1}}}\|\D_{l}a^{n}\|_{L^{\infty}_{T}(L^{p_{1}})}\leq C\big(
\sup_{l\geq m}2^{l\frac{N}{p_{1}}}\|\D_{l}a_{0}\|_{L^{p_{1}}})\big)\big(\|\tilde{u}^{n}\|_{\widetilde{L}^{1}_{t}(B^{\frac{N}{p_{2}}}_{p_{2},r})}+
\|u^{n}_{L}\|_{\widetilde{L}^{1}_{t}(B^{\frac{N}{p_{2}}}_{p_{2},r})}\big)
\label{39}
\end{equation}
Using (\ref{1conduti}) and $({\cal H}_{5})$, $({\cal H}_{6})$, we thus get:
$$\|a^{n}-S_{m}a^{n}\|_{L^{\infty}_{T}(B^{\frac{N}{p_{1}}}_{p_{1},\infty})\cap L^{\infty}}
\leq\sup_{l\geq m}2^{l\frac{N}{p_{1}}}\|\D_{l}a_{0}\|_{L^{p_{1}}}+\frac{C}{\log2}
(1+\|a_{0}\|_{B^{\frac{N}{p_{1}}}_{p_{1},\infty}\cap L^{\infty}})(1+\underline{\nu}^{-1}\widetilde{L}_{0})\eta.$$
Hence $({\cal H}_{1})$ is strictly satisfied provided that $\eta$ further satisfies:
\begin{equation}
\frac{C}{\log2}(1+\|a_{0}\|_{B^{\frac{N}{p_{1}}}_{p_{1},\infty}\cap L^{\infty}})(1+\underline{\nu}^{-1}\widetilde{L}_{0})\eta<\frac{c\underline{\nu}}{2\bar{\nu}}.
\label{40}
\end{equation}
So ${\cal H}_{1}$ is strictly verified.\\
$({\cal H}_{3})$ is trivially verified by the transport equation as we assume that $1+a_{0}$ is bounded away and that $a_{0}\in L^{\infty}$.
\\
Next, according to $({\cal H}_{4})$ condition $({\cal H}_{2})$ is satisfied provide:
\begin{equation}
 T<\frac{2^{-2m}\underline{\nu}}{C\bar{\nu}^{2}A_{0}^{2}}
\label{41}
\end{equation}
One can now conclude that if $T<T^{n}$ has been choosen so that conditions (\ref{36}), (\ref{38}) and (\ref{41})  are satisfied (with $\eta$ verifying (\ref{1conduti}) and (\ref{40}), and $m$ defined in (\ref{def})
and $n\geq m$ then $(a^{n},u^{n},\Pi^{n})$ satisfies $({\cal H}_{1})$ to $({\cal H}_{8})$, thus is bounded independently of $n$
on $[0,T]$.\\
We still have to state that $T^{n}$ may be bounded by below by the supremum $\bar{T}$ of all times $T$
such that (\ref{36}), (\ref{38}) and (\ref{41})  are satisfied. This is actually a consequence of the uniform bounds we have just obtained, and of a theorem of blow-up of R. Danchin in \cite{Da4}. Indeed, by combining all these informations, one can prove that if $T^{n}<\bar{T}$ then $(a^{n},u^{n},\n\Pi^{n})$ is actually in:
$$
\begin{aligned}
&\widetilde{L}^{\infty}_{T^{n}}(B^{\frac{N}{p_{1}}+\e^{'}}_{p_{1},\infty}\cap B^{\N+1}_{2,1})\times\big(\widetilde{L}^{\infty}_{T^{n}}(B^{\frac{N}{p_{2}}-1}_{p_{2},r}\cap B^{\N}_{2,1})\cap \widetilde{L}^{1}_{T^{n}}(B^{\frac{N}{p_{2}}+1}_{p_{2},r}\cap B^{\N+2}_{2,1})\big)^{N}\\
&\hspace{9cm}\times\widetilde{L}^{1}_{T^{n}}(B^{\frac{N}{p_{2}}-1}_{p_{2},r}\cap B^{\N}_{2,1}).
\end{aligned}
$$
hence may be continued beyond $\bar{T}$ (see the remark on the lifespan
in \cite{Da4} where a control of $\n u$ in $\widetilde{L}^{1}(B^{0}_{\infty,\infty})$ is required). We thus have $T^{n}\geq\bar{T}$.
\subsubsection*{Compactness arguments}
We now have to prove that $(a^{n},u^{n})_{n\in\mathbb{N}}$ tends (up to a subsequence) to some
function $(a,u)$ which belongs to $E_{T}$ and satisfies $(\ref{0.6})$.
The proof is based on Ascoli's theorem and compact embedding for Besov spaces.
As similar arguments have been employed in \cite{Da2} or \cite{Da3}, we only give the outlines of the proof.
\begin{itemize}
\item Convergence of $(a^{n})_{n\in\mathbb{N}}$:\\
We use the fact that $\widetilde{a}^{n}=a^{n}-a^{n}_{0}$ satisfies:
$$\p_{t}\widetilde{a}^{n}=-u^{n}\cdot\n a^{n}.$$
Since $(u^{n})_{n\in\mathbb{N}}$ is uniformly bounded in $\widetilde{L}^{1}_{T}(B^{\frac{N}{p_{2}}+1}_{p_{2},r})
\cap \widetilde{L}^{\infty}_{T}(B^{\frac{N}{p_{2}}-1}_{p_{2},r})$, it is, by interpolation,
also bounded in $\widetilde{L}^{r^{'}}_{T}(B^{\frac{N}{p_{2}}-1+\frac{2}{r^{'}}}_{p_{2},r})$
for any $r^{'}\in[1,+\infty]$.
By taking $r=2$ and using the standard product laws in Besov spaces, we thus easily gather that
$(\p_{t}\widetilde{a}^{n})$ is uniformly bounded in $\widetilde{L}^{2}_{T}(B^{\frac{N}{p_{1}}-1}_{p_{1},\infty})$.
$$
\begin{aligned}
\|\p_{t}\widetilde{a}^{n}\|_{\widetilde{L}^{2}_{T}(B^{\frac{N}{p_{1}}-1}_{p_{1},\infty})}&\lesssim
\|u^{n}\|_{\widetilde{L}^{2}_{T}(B^{\frac{N}{p_{2}}}_{p_{2},r})}\|\n a^{n}\|_{\widetilde{L}^{\infty}_{T}(B^{
\frac{N}{p_{1}}-1}_{p_{1},\infty})}.
\end{aligned}
$$
Hence $(\widetilde{a}^{n})_{n\in\mathbb{N}}$ is bounded in
$\widetilde{L}^{\infty}_{T}(B^{\frac{N}{p_{1}}-1}_{p_{1},\infty}\cap B^{\frac{N}{p_{1}}}_{p_{1},\infty})$
and equicontinuous on $[0,T]$ with values in $B^{\frac{N}{p_{1}}-1}_{p_{1},\infty}$. Since the embedding
$B^{\frac{N}{p_{1}}-1}_{p_{1},\infty}\cap B^{\frac{N}{p_{1}}}_{p_{1},\infty}\hookrightarrow
B^{\frac{N}{p_{1}}-1}_{p_{1},\infty}$ is (locally) compact, and $(a_{0}^{n})_{n\in\mathbb{N}}$ tends
to $a_{0}$ in $B^{\frac{N}{p_{1}}}_{p_{1},\infty}$, we conclude that $(a^{n})_{n\in\mathbb{N}}$ tends (up to extraction)
to some distribution $a$. Given that $(a^{n})_{n\in\mathbb{N}}$ is bounded in
$\widetilde{L}^{\infty}_{T}(B^{\frac{N}{p_{1}}+\frac{\e}{2}}_{p_{1},r})$,
we actually have $a\in\widetilde{L}^{\infty}_{T}(B^{\frac{N}{p_{1}}+\frac{\e}{2}}_{p_{1},r})$.
\item Convergence of $(u^{n}_{L})_{n\in\mathbb{N}}$:\\
From the definition of $u^{n}_{L}$ and classical proposition on Stokes equation, it is clear that $(u^{n}_{L})_{n\in\mathbb{N}}$ and
$(\n\Pi^{n}_{L})_{n\in\mathbb{N}}$
tend to solution $u_{L}$ and $\n\Pi_{L}$ to:
$$\p_{t}u_{L}-\mu\D u_{l}+\n\Pi_{L}=f,\;\;u_{L}(0)=u_{0}$$
in $\widetilde{L}^{\infty}_{T}(B^{\frac{N}{p_{2}}-1}_{p_{2},r})\cap \widetilde{L}^{1}_{T}(B^{\frac{N}{p_{2}}+1}
_{p_{2},r})$ for $(u^{n}_{L})_{n\in\mathbb{N}}$ and $\widetilde{L}^{1}_{T}(B^{\frac{N}{p_{2}}-1}
_{p_{2},r})$ for $(\n\Pi^{n}_{L})_{n\in\mathbb{N}}$.
\item  Convergence of $(\widetilde{u}^{n})_{n\in\mathbb{N}}$:\\
We use the fact that:
$$\p_{t}\widetilde{u}^{n}=-u^{n}_{L}\cdot\n\widetilde{u}^{n}-\widetilde{u}^{n}\cdot\n u^{n}+
(1+a^{n})\D\widetilde{u}^{n}+a^{n}\D u^{n}_{L}-u^{n}_{L}\cdot\n u^{n}_{L}-\n\widetilde{\Pi}^{n}.$$
As $(a^{n})_{n\in\mathbb{N}}$ is uniformly bounded in $\widetilde{L}^{\infty}_{T}(B^{\frac{N}{p_{1}}}_{p_{1},\infty})$
and $(u^{n})_{n\in\mathbb{N}}$ is uniformly bounded in $\widetilde{L}^{\infty}_{T}(B^{\frac{N}{p_{2}}-1}_{p_{2},r})
\cap \widetilde{L}^{\frac{4}{3}}(B^{\frac{N}{p_{2}}+\frac{1}{2}}_{p_{2},r})$, it is easy to see that the last term
of the right-hand side is uniformly bounded in $\widetilde{L}^{\infty}_{T}(B^{\frac{N}{p_{2}}-1}_{p_{2},r})$ and that
the other terms are
uniformly bounded in $\widetilde{L}^{\frac{4}{3}}(B^{\frac{N}{p_{2}}-\frac{3}{2}}_{p_{2},r})$.\\
Hence $(\widetilde{u}^{n})_{n\in\mathbb{N}}$ is bounded in $\widetilde{L}^{\infty}_{T}(B^{\frac{N}{p_{2}}-1}_{p_{2},r})$
and equicontinuous on $[0,T]$ with values in $B^{\frac{N}{p_{2}}-1}_{p_{2},r}+B^{\frac{N}{p_{2}}-\frac{3}{2}}_{p_{2},r}$.
This enables to conclude that
$(\widetilde{u}^{n})_{n\in\mathbb{N}}$ converges (up to extraction) to some function $\widetilde{u}\in
\widetilde{L}^{\infty}_{T}(B^{\frac{N}{p_{2}}-1}_{p_{2},r})\cap \widetilde{L}^{1}_{T}(B^{\frac{N}{p_{2}}+1}_{p_{2},r})$.
\end{itemize}
We proceed similarly for $(\Pi^{n}_{L})_{n\in\mathbb{N}}$ and $(\widetilde{\Pi}^{n})_{n\in\mathbb{N}}$.
By interpolating with the bounds provided by the previous step, one obtains better
results of convergence so that one can pass to the limit in the mass equation and in (\ref{3momentum}).
Finally by setting $u=\widetilde{u}+u_{L}$ and $\Pi=\widetilde{\Pi}+\Pi_{L}$, we conclude that
$(a,u,\Pi)$ satisfies (\ref{0.6}).\\
In order to prove continuity in time for $a$ it suffices to make use of proposition
\ref{transport1}. Indeed, $a_{0}$ is in $B^{\frac{N}{p_{1}}+\e}_{p_{1},\infty}\cap L^{\infty}$, and having $a\in \widetilde{L}^{\infty}_{T}(B^{\frac{N}{p_{1}}+\frac{\e}{2}}
_{p_{1},\infty})\cap L^{\infty}$
and $u\in \widetilde{L}^{1}_{T}(B^{\frac{N}{p_{2}}+1}_{p_{2},r})$ insure that
$\p_{t}a+u\cdot\n a$ belongs to $\widetilde{L}^{1}_{T}(B^{\frac{N}{p_{1}}}_{p_{1},\infty})$.
Similarly, continuity for $u$ may be proved by using that $u_{0}\in B^{\frac{N}{p_{2}}}_{p_{2},r}$
and that $(\p_{t}u-\mu\D u)\in \widetilde{L}^{1}_{T}(B^{\frac{N}{p_{2}}-1}_{p_{2},r})$.
\section{The proof of the uniqueness}
\label{section6}
\subsection{Uniqueness when $1\leq p_{2}<2N$, $\frac{2}{N}<\frac{1}{p_{1}}+\frac{1}{p_{2}}$ and $N\geq3$.}
In this section, we focus on the case $N\geq3$ and postpone the analysis of the other cases
(which turns out to be critical) to the next section.
Throughout the proof, we assume that we are given two solutions $(a^{1},u^{1})$ and $(a^{2},u^{2})$ of (\ref{0.6})
which belongs to:
$$ \big(\widetilde{C}([0,T]; B^{\frac{N}{p_{1}}+\frac{\e}{2}}_{p_{1},\infty})\cap L^{\infty}\big)\times\big(\widetilde{C}([0,T]; B^
{\frac{N}{p_{2}}-1}_{p_{2},r})
\cap \widetilde{L}^{1}(0,T;B^{\frac{N}{p_{2}}+1}_{p_{2},r})\big)^{N}.$$
Let $\delta a=a^{2}-a^{1}$, $\delta u=u^{2}-u^{1}$ and $\delta\Pi=\Pi^{2}-\Pi^{1}$. The system for $(\delta a,\delta u)$ reads:
\begin{equation}
\begin{cases}
\begin{aligned}
&\p_{t}\delta a+u^{2}\cdot\n\delta a=-\delta u\cdot\n a^{1},\\
&\p_{t}\delta u-(1+a^{2})(\mu\D\delta u-\n\delta\Pi)=
F(a^{i},u^{i},\Pi^{i}).\\
\end{aligned}
\end{cases}
\label{43}
\end{equation}
with:
$$F(a^{i},u^{i},\Pi^{i})=u^{1}\cdot\n\delta u+\delta u\cdot\n u^{2}+\delta a(\mu\D u^{1}-\n\Pi^{1}).$$
The function $\delta a$ may be estimated by taking advantage of proposition \ref{transport1} with $s=\frac{N}{p_{1}}-1+\frac{\e}{2}$.
We get for all
$t\in[0,T]$,
$$\|\delta a(t)\|_{B^{\frac{N}{p_{1}}-1}_{p_{1},\infty}}\leq
C\|\delta u\cdot\n a^{1}\|_{\widetilde{L}_{T}^{1}(B^{\frac{N}{p_{1}}-1+\frac{\e}{2}}_{p_{1},\infty})}\exp\big(\frac{C}{\e^{r-1}}
(V_{p_{1},1-\frac{1}{r}}(t))^{r}\big),$$
We have then by proposition \ref{produit1} and \ref{propChemin2}:
\begin{equation}
\|\delta a(t)\|_{B^{\frac{N}{p_{1}}-1}_{p_{1},\infty}}\leq
C\|\delta u\|_{\widetilde{L}_{T}^{1}(B^{\frac{N}{p_{2}}}_{p_{2},r})}\|a^{1}\|_{\widetilde{L}_{T}^{\infty}(B^{\frac{N}{p_{1}}+\frac{\e}{2}}_
{p_{1},\infty})}\exp\big(\frac{C}{\e^{r-1}}
(\|u^{2}\|_{\widetilde{L}_{T}^{1}(B^{\frac{N}{p_{2}}+1}_{p_{2},r})})^{r}\big),
\end{equation}
For bounding $\delta u$, we aim at applying proposition \ref{linearise} to the second equation of (\ref{43}).
So let us fix an integer $m$ such that:
\begin{equation}
1+\inf_{(t,x)\in[0,T]\times\R^{N}}S_{m}a^{2}\geq\frac{\underline{b}}{2}\;\;\mbox{and}\;\;\|a^{2}-
S_{m}a^{2}\|_{L^{\infty}_{T}(B^{\frac{N}{p_{1}}}_{p_{1},\infty})}\leq c\frac{\underline{\nu}}{\bar{\nu}}.
\label{45}
\end{equation}
Now applying proposition \ref{linearise} with $s=\frac{N}{p_{2}}-2$ insures that for all time $t\in[0,T]$, we have:
\begin{equation}
\begin{aligned}
&\|u\|_{\widetilde{L}^{\infty}_{T}(B^{\frac{N}{p_{2}}-2}_{p_{2},r})}+\kappa\underline{\nu}
\|u\|_{\widetilde{L}^{1}_{T}(B^{\frac{N}{p_{2}}}_{p_{2},r})}+\|\n\delta\Pi\|_{\widetilde{L}^{1}_{T}(B^{\frac{N}{p_{2}}-2}_{p_{2},r})}\leq e^{CZ_{m}(T)}\times\\
&\hspace{3cm}\big(\|{\cal P}F(a^{i},u^{i},\Pi^{i})\|_{\widetilde{L}^{1}_{T}(B^{\frac{N}{p_{2}}-2}_{p_{2},r})}+(\frac{\underline{\nu}(p_{2}-1)}{p_{2}}){\cal A}_{T}\|u\|_{\widetilde{L}^{1}_{T}
(B^{\frac{N}{p_{2}}-\alpha^{'}}_{p,r})}\big).\\
\end{aligned}
\label{3impo1}
\end{equation}
with $Z_{m}(t)=2^{m}\mu^{2}\underline{\nu}^{-1}\int^{t}_{0}\|a(\tau)\|^{2}_{B^{\frac{N}{p_{1}}}_{p_{1},\infty}\cap L^{\infty}}d\tau$.\\
Hence, applying proposition \ref{produit1}, corollary \ref{produit2} and the fact that ${\rm div}\delta u=0$, we get as exemple:
$$
\begin{aligned}
&\|\delta u\cdot\n u^{1}\|_{\widetilde{L}^{1}_{T}(B^{\frac{N}{p_{2}}-2}_{p_{2},r})}\lesssim \|u^{1}\|_{\widetilde{L}^{\infty}_{T}(B^{\frac{N}{p_{2}}-1}_{p_{2},r})}^{\frac{1}{2}}
\|u^{1}\|_{\widetilde{L}^{1}_{T}(B^{\frac{N}{p_{2}}+1}_{p_{2},r})}^{\frac{1}{2}}\big(
\|\delta u\|_{\widetilde{L}^{\infty}_{T}(B^{\frac{N}{p_{2}}-2}_{p_{2},r})}+\|\delta u\|_{\widetilde{L}^{\infty}_{T}(B^{\frac{N}{p_{2}}}_{p_{2},r})}\big).
\end{aligned}
$$By the fact that $\frac{2}{N}<\frac{1}{p_{1}}+\frac{1}{p_{2}}$, $N\geq 3$ and $\frac{1}{p_{2}}\leq\frac{1}{N}+\frac{1}{p_{1}}$ imply that:
\begin{equation}
\begin{aligned}
&\|\delta a(\mu\D u_{1}-\n\delta\Pi)\|_{\widetilde{L}^{1}_{T}(B^{\frac{N}{p_{2}}-2}_{p_{2},r})}\lesssim
\|\delta a\|_{\widetilde{L}^{\infty}_{T}(B^{\frac{N}{p_{1}}-1}_{p_{1},\infty})}\big(\|\D u_{1}\|_{\widetilde{L}^{1}_{T}(B^{\frac{N}{p_{2}}-1}_{p_{1},r})}\\
&\hspace{9cm}+\|\n\Pi^{1}\|_{\widetilde{L}^{1}_{T}(B^{\frac{N}{p_{2}}-1}_{p_{2},
r})}\big).
\end{aligned}
\label{3impo}
\end{equation}
Now let choose $T_{1}$ enough small to controll in (\ref{3impo1}) ${\cal A}_{T}\|u\|_{\widetilde{L}^{1}_{T}
(B^{\frac{N}{p_{2}}-\alpha^{'}}_{p,r})}$ and in (\ref{3impo}) $\|\D u_{1}\|_{\widetilde{L}^{1}_{T}(B^{\frac{N}{p_{2}}-1}_{p_{1},r})}+\|\n\Pi^{1}\|_{\widetilde{L}^{1}_{T}(B^{\frac{N}{p_{2}}-1}_{p_{2},
r})}$ and by the fact that $\|a_{1}\|_{\widetilde{L}^{\infty}_{T}(B^{\frac{N}{p_{1}}+\frac{\e}{2}}_{p_{1},\infty})}\leq c$ with $c$ small, we obtain finally:
$$\|(\delta a,\delta u,\n\delta\Pi)\|_{F_{T_{1}}}\leq c C\|(\delta a,\delta u,\n\delta\Pi)\|_{F_{T_{1}}},$$
with:
$$F_{T}=\widetilde{C}_{T}([0,T],B^{\frac{N}{p_{1}}}_{p_{1},\infty})\times\big(\widetilde{L}^{1}_{T}(B^{\frac{N}{p_{2}}}_{p_{2},r})
\cap \widetilde{L}^{\infty}_{T}(B^{\frac{N}{p_{2}}-2}_{p_{2},r})\big)\times \widetilde{L}^{1}_{T}(B^{\frac{N}{p_{2}}-2}_{p_{2},r}).$$
We obtain so $(\delta a,\delta u,\n \delta\Pi)=0$ on $[0,T_{1}]$ for $T_{1}$ enough small. By connectivity we obtain that
$(\delta a,\delta u,\delta\n\Pi)=0$ on $[0,T]$. This conclude this case.
\subsection{Uniqueness when:$\frac{2}{N}=\frac{1}{p_{1}}+\frac{1}{p_{2}}$ or $p_{2}=2N$ or $N=2$.}
The above proof fails in dimension two. One of the reasons why is that the product of functions does not map
$B^{\frac{N}{p_{1}}-1}_{p_{1},\infty}\times B^{\frac{N}{p_{2}}-1}_{p_{2},r}$ in $B^{\frac{N}{p_{2}}-2}_{p_{2},r}$ but only in the larger space $B^{-1}_{2,\infty}$.
This induces us to bound $\delta a$ in $\widetilde{L}_{T}^{\infty}(B^{\frac{N}{p_{1}}-1}_{p_{1},\infty})$ and $\delta u$ in
$\widetilde{L}_{T}^{\infty}(B^{\frac{N}{p_{2}}-2}_{p_{2},\infty})\cap \widetilde{L}_{T}^{1}(B^{\frac{N}{p_{2}}}_{p_{2},\infty})$.
In fact, it is enough to study only the case $\frac{2}{N}=\frac{1}{p_{1}}+\frac{1}{p_{2}}$. Indeed the other cases deduct from this case.
If $p_{2}=2N$ then $p_{1}=\frac{2N}{3}$ as $\frac{1}{p_{1}}\leq\frac{1}{N}+\frac{1}{p_{2}}$ and $\frac{2}{N}\leq \frac{1}{p_{1}}+\frac{1}{p_{2}}$. So it is a particular case of $\frac{2}{N}=\frac{1}{p_{1}}+\frac{1}{p_{2}}$. For $N=2$, we begin with $p_{2}=4$ and $p_{1}=\frac{4}{3}$ and by embedding we get the result for $1\leq p_{1}\leq\frac{4}{3}$,  $1\leq p_{2}\leq 4$
and for $1\leq p_{1}\leq 4$,  $1\leq p_{2}\leq \frac{4}{3}$.\\
Moreover in your case, it exists two possibilities, one when $1<p_{2}<2$ and when $p_{2}\geq2$. The first case is resolved by embedding so we have just to treat the case $p_{2}\geq 2$. We want show that $(\delta a,\delta u,\n\delta\Pi)\in G_{T}$ where:
$$G_{T}=\widetilde{C}_{T}([0,T],B^{\frac{N}{p_{1}}-1+\frac{\e}{2}}_{p_{1},\infty})\times\big(\widetilde{L}^{1}_{T}(B^{\frac{N}{p_{2}}}_{p_{2},\infty})
\cap \widetilde{L}^{\infty}_{T}(B^{\frac{N}{p_{2}}-2}_{p_{2},\infty})\big)\times \widetilde{L}^{1}_{T}(B^{\frac{N}{p_{2}}-2}_{p_{2},\infty}).$$
In fact we proceed exactly as in the previous proof but we benefit that $a$ is in $\widetilde{C}_{T}([0,T],B^{\frac{N}{p_{1}}-1+\frac{\e}{2}}_{p_{1},\infty})$ which give sense to the product $\delta\D u$. It conclude the proof of uniqueness.
\subsection{Proof of theorem \ref{theo11}}
The proof follow strictly the same lines than the proof of theorem \ref{theo1} except that we profits of the fact
that in the proposition \ref{linearise} in the case of $p_{2}=2$ we do not need of conditions of smallness on the initial density $a_{0}$.
\section{Proof of theorem \ref{theo2} and \ref{theo22}}
In this case we have a control on $u$ only in $\widetilde{L}^{1}(B^{\frac{N}{p_{2}}+1}_{p_{2},\infty})$, that is why to control the density in this case we need to use the proposition \ref{transport2}. The rest of the proofs is similar to proof \ref{theo1}.
\label{section7}
\section{Appendix}
\label{section9}
\subsection{Elliptic estimates}
This section is devoted to the study of the elliptic equation:
\begin{equation}
{\rm div}(b\n\Pi)={\rm div}F.
\label{ellipticappendice}
\end{equation}
with $b=1+a$.\\
Let us first study the stationary case where $F$ and $b$ are independent of the time:
\begin{proposition}
\label{propositionA5}
Let $0<\alpha<1$, $(p,r)\in[1,+\infty]^{2}$ and $\sigma\in\R$ satisfy $\alpha\leq\sigma\leq\alpha+\frac{N}{p_{1}}$. Then the operator
${\cal H}_{b}:F\rightarrow\n\Pi$ is a linear bounded operator in $B^{\sigma}_{p,r}$ and the following estimate holds true:
\begin{equation}
\underline{b}\|\n\Pi\|_{B^{\sigma}_{p,r}}\lesssim{\cal A}^{\frac{|\sigma|}{\min(1,\alpha)}}\|{\cal Q}F\|_{B^{\sigma}_{p,r}},
\label{63a}
\end{equation}
with if $\alpha\ne 1$:
$${\cal A}=1+\underline{b}^{-1}\|\n b\|_{B^{\frac{N}{p_{1}}+\alpha-1}_{p_{1},r}}.$$
\end{proposition}
{\bf Proof:}\\
\\
Let us first rewrite (\ref{ellipticappendice}) as follows:
\begin{equation}
{\rm div}(b_{m}\n\Pi)={\rm div}F-E_{m},
\label{ellipticappendice1}
\end{equation}
with $E_{m}={\rm div}((Id-S_{m})a\n\Pi)$.\\
Apply $\D_{q}$ to (\ref{ellipticappendice1}) we get:
\begin{equation}
{\rm div}(b_{m}\n\D_{q}\Pi)={\rm div}F_{q}-\D_{q}(E_{m})+R_{q},
\label{ellipticappendice2}
\end{equation}
with $R_{q}={\rm div}(b_{m}\n\Pi_{q})-\D_{q}{\rm div}(b_{m}\n\Pi)$.
Multiplying (\ref{ellipticappendice2}) by $\D_{q}\Pi|\D_{q}\Pi|^{p-2}$ and integrate, we gather:
\begin{equation}
\begin{aligned}
&\int_{\R^{N}}b_{m}|\n\Pi_{q}|^{2}|\Pi_{q}|^{p-2}dx+\int_{\R^{N}}b_{m}|\n |\Pi_{q}|^{2}|^{p}dx
\leq(\|{\rm div}F_{q}\|_{L^{p}}+\|R_{q}\|_{L^{p}})\|\Pi_{q}\|_{L^{p}}^{p-1}\\
&+\int_{\R^{N}}|(Id-S_{m})a|\,|\n\Pi_{q}|^{2}|\Pi_{q}|^{p-2}dx+\int_{\R^{N}}|(Id-S_{m})a|\,|\n |\Pi_{q}|^{2}|^{p}dx,
\label{64a}
\end{aligned}
\end{equation}
Assuming that $m$ has been choose so large as to satisfy:
$$\|a-S_{m}a\|_{L^{\infty}_{T}(B^{\frac{N}{p_{1}}+\alpha}_{p_{1},\infty})\cap L^{\infty}}\leq\frac{\underline{b}}{2},$$
an by using lemma A5 in \cite{Da6}:
$$2^{2q}\|\Pi_{q}\|_{L^{p}}\lesssim 2^{q}\|F_{q}\|_{L^{p}}+\|R_{q}\|_{L^{p}}.$$
By multipliyng by $2^{q(s-1)}$ and by integrating on $l^{r}$ we get:
$$\underline{b}\|\n\Pi\|_{B^{s}_{p,r}}\lesssim\|{\cal Q}f\|_{B^{s}_{p,r}}+\|R_{q}\|_{B^{s}_{p,r}}.$$
The commutator may be bounded thanks to lemma \ref{alemme3} with $0<\alpha<1$, $\sigma=s-1$. We have get:
$$\underline{b}\|\n\Pi\|_{B^{s}_{p,r}}\lesssim\|{\cal Q}F\|_{B^{s}_{p,r}}+\|\n S_{m}a\|
_{B^{\frac{N}{p_{1}}+\alpha}_{p_{1},r}}
\|\n\Pi\|_{B^{s-\alpha}_{p,r}}.$$
Therefore complex interpolation entails:
$$\|\n\Pi\|_{B^{s-\alpha}_{p,r}}\leq\|\n\Pi\|_{B^{s}_{p,r}}^{\frac{s-\alpha}{s}}
\|\n\Pi\|^{\frac{\alpha}{s}}_{B^{0}_{p,\infty}}.$$ 
Note that, owing to Bersntein inequality, we have:
$$\|\n S_{m}a\|_{B^{\frac{N}{p_{1}}+\alpha}_{p_{1},r}}\lesssim
2^{m+\alpha}\|a\|_{B^{\frac{N}{p_{1}}}_{p_{1},r}}.$$
We have then:
$$\underline{b}\|\n\Pi\|_{B^{s}_{p,r}}\lesssim\|{\cal Q}f\|_{B^{s}_{p,r}}+
2^{m+\alpha}\|a\|_{B^{\frac{N}{p_{1}}}_{p_{1},\infty}}
\|\n\Pi\|_{B^{s}_{p,r}}^{\frac{s-\alpha}{s}}\|\n\Pi\|^{\frac{\alpha}{s}}_{B^{0}_{p,\infty}}.$$
And we conclude by Young's inequality with $p_{1}=\frac{s}{s-\alpha}$ and  $p_{2}=\frac{s}{\alpha}$. And we recall
\subsection{Commutator estimates}
This section is devoted to the proof of commutator estimates which have been used in section $2$ and $3$.
They are based on
paradifferentiel calculus, a tool introduced by J.-M. Bony in \cite{Bo}. The basic idea of
paradifferential calculus is that
any product of two distributions $u$ and $v$ can be formally decomposed into:
$$uv=T_{u}v+T_{v}u+R(u,v)=T_{u}v+T^{'}_{v}u$$
where the paraproduct operator is defined by $T_{u}v=\sum_{q}S_{q-1}u\D_{q}v$, the remainder operator $R$, by
$R(u,v)=\sum_{q}\D_{q}u(\D_{q-1}v+\D_{q}v+\D_{q+1}v)$ and $T^{'}_{v}u=T_{v}u+R(u,v)$.
\\
Inequality (\ref{13}) is a consequence of the following lemma:
\begin{lemme}
\label{alemme3}
Let $p_{1}\in[1,+\infty]$, $p\in[1,+\infty]$,
$\alpha\in(1-\NN,1[$, $k\in\{1,\cdots,N\}$ and $R_{q}=\D_{q}(a\p_{k}w)-\p_{k}(a\D_{q}w)$. There
exists $c=c(\alpha,N,\sigma)$ such that:
\begin{equation}
2^{q\sigma}\|\widetilde{R}_{q}\|_{L^{p}}\leq C c_{q}\|a\|_{B^{\frac{N}{p_{1}}+\alpha}_{p_{1},r}}\|w\|_{B^{\sigma+1-\alpha}_{p,r}}
\label{57}
\end{equation}
whenever $-\frac{N}{p_{1}}<\sigma\leq \frac{N}{p_{1}}+\alpha$ and where $c_{q}\in L^{r}$.\\
In the limit case $\sigma=-\frac{N}{p_{1}}$, we have for some constant $C=C(\alpha,N)$:
\begin{equation}
2^{-q\frac{N}{p_{1}}}\|\widetilde{R}_{q}\|_{L^{p}}\leq C\|a\|_{B^{\alpha+\frac{N}{p_{1}}}_{p,1}}\|w\|_{B^{-\frac{N}{p_{1}}+1-\alpha}_{p,\infty}}.
\label{58}
\end{equation}
\end{lemme}
{\bf Proof}\\
\\
The proof is based on Bony's decomposition which enables us to split $R_{q}$
into:
$$R_{q}=\underbrace{\p_{k}[\D_{q},T_{a}]w}_{R_{q}^{1}}-\underbrace{\D_{q}T_{\p_{k}a}w}_{R_{q}^{2}}+\underbrace{\D_{q}T_{\p_{k}w}a}_{R_{q}^{3}}
+\underbrace{\D_{q}R(\p_{k}w,a)}_{R_{q}^{4}}-\underbrace{\p_{k}T^{'}_{\D_{q}w}a}_{R_{q}^{5}}.$$
By using the fact that:
$$R^{1}_{q}=\sum^{q+4}_{q^{'}=q-4}\p_{k}[\D_{q},S_{q^{'}-1}a]\D_{q^{'}}w,$$
Using the definition of the operator $\D_{q}$ leads to:
$$[\D_{q},S_{q^{'}-1}a]\D_{q^{'}}w(x)=-\int h(y)\big(S_{q^{'}-1}a(x)-
S_{q^{'}-1}a(x-2^{-q}y)\big)\D_{q^{'}}w(x-2^{-q}y)dy.$$
and:
$$
\begin{aligned}
\big|[\D_{q},S_{q^{'}-1}a]\D_{q^{'}}w(x)\big|&\leq \|\n S_{q^{'}-1}a\|_{L^{\infty}}2^{-q}\int 2^{qN}|h(2^{q}u)||2^{q}u|
|\D_{q^{'}}w|(x-u)du,\\
&\leq 2^{qN}\|\n S_{q^{'}-1}a\|_{L^{\infty}}|\big(h(2^{q}\cdot)|\cdot|*\D_{q^{'}}w\big)|(x).
\end{aligned}
$$
So we get:
$$\|[\D_{q},S_{q^{'}-1}a]\D_{q^{'}}w\|_{L^{p}}\leq\|\n S_{q^{'}-1}a\|_{L^{\infty}}\|\D_{q^{'}}w\|_{L^{p}}$$
we readily get under the hypothesis that $\alpha<1$,
\begin{equation}
2^{q\sigma}\|R^{1}_{q}\|_{L^{p}}\lesssim\sum^{q+4}_{q^{'}=q-4}2^{q\sigma}
\|\n S_{q^{'}-1}a\|_{L^{\infty}}\|\D_{q^{'}}w\|_{L^{p}}.
\label{59bis}
\end{equation}
We have then:
\begin{equation}
2^{q\sigma}\|R^{1}_{q}\|_{L^{p}}\lesssim c_{q}\|\n a\|_{B^{\alpha-1}_{\infty,\infty}}\|w\|_{B^{\sigma+1-\alpha}_{p,r}}.
\label{59}
\end{equation}
In the case $\alpha=1$, we get:
\begin{equation}
2^{q\sigma}\|R^{1}_{q}\|_{L^{p}}\lesssim c_{q}\|\n a\|_{B^{0}_{\infty,1}}\|w\|_{B^{\sigma+1-\alpha}_{p,r}}.
\label{591}
\end{equation}
For bounding $R^{2}_{q}$, standard continuity results for the paraproduct insure that if $\alpha<1$, $R^{2}_{q}$
satisfies that:
$$2^{q\sigma}\|R^{2}_{q}\|_{L^{p}}\leq c_{q}\|\n a\|_{B^{\alpha-1}_{\infty,\infty}}\|w\|_{B^{\sigma+1-\alpha}_{p,r}}.$$
and if $\alpha=1$
$$2^{q\sigma}\|R^{2}_{q}\|_{L^{p}}\leq c_{q}\|\n a\|_{B^{\alpha-1}_{\infty,1}}\|w\|_{B^{\sigma+1-\alpha}_{p,r}}.$$
Standard continuity results for the paraproduct insure that $R^{3}_{q}$ satisfies:
\begin{equation}
2^{q\sigma}\|R^{3}_{q}\|_{L^{p}}\lesssim c_{q}\|\n w\|_{B^{\sigma-\alpha-\NN}_{\infty,\infty}}
\|a\|_{B^{\NN+\alpha}_{p,r}}.
\label{60}
\end{equation}
provided $\sigma-\alpha-\NN<0.$\\
If $\sigma-\alpha-\NN=0$ then:
\begin{equation}
2^{q\sigma}\|R^{3}_{q}\|_{L^{p}}\lesssim c_{q}\|\n w\|_{B^{0}_{\infty,1}}
\|a\|_{B^{\NN+\alpha}_{p,r}}.
\label{601}
\end{equation}
Next, standard continuity result for the remainder insure that under the hypothesis $\sigma>-\NN$, we have:
\begin{equation}
2^{q\sigma}\|R^{4}_{q}\|_{L^{p}}\lesssim c_{q}\|\n w\|_{B^{\sigma-\alpha}_{p,r}}\|a\|_{B^{\NN+\alpha}_{p,\infty}}.
\label{61}
\end{equation}
For bounding $R^{5}_{q}$ we use the decomposition:
$$R^{5}_{q}=\sum_{q^{'}\geq q-3}\p_{k}(S_{q^{'}+2}\D_{q}w\D_{q^{'}}a),$$
which leads (after a suitable use of Bernstein and H\"older inequalities) to:
$$
\begin{aligned}
2^{q\sigma}\|R^{5}_{q}\|_{L^{p}}&\lesssim\sum_{q^{'}\geq q-3}2^{q\sigma}2^{q^{'}}\|\D_{q^{'}}a\|_{L^{\infty}}
\|S_{q^{'}+2}\D_{q}w\|_{L^{p}}\\
&\lesssim\sum_{q^{'}\geq q-2}2^{(q-q^{'})(\alpha+\NN-1)}2^{q(\sigma+1-\alpha)}
\|\D_{q}w\|_{L^{p}}2^{q^{'}(\NN+\alpha)}
\|\D_{q^{'}}a\|_{L^{p}}.
\end{aligned}
$$
Hence, since $\alpha+\NN-1>0$, we have:
$$2^{q\sigma}\|R^{5}_{q}\|_{L^{p}}\lesssim c_{q}\|\n w\|_{B^{\sigma+1-\alpha}_{p,r}}\|a\|_{B^{\NN+\alpha}_{p,\infty}}.$$
Combining this latter inequality with (\ref{59}), (\ref{60}) and (\ref{61}), and using the embedding
$B^{\NN}_{p,r}\hookrightarrow B^{r-\NN}_{\infty,\infty}$
for $r=\NN+\alpha-1$, $\sigma-\alpha$ completes the proof of (\ref{57}).\\
\\
\\
The proof of (\ref{58}) is almost the same: for bounding $R^{1}_{q}$, $R^{2}_{q}$, $R^{3}_{q}$ and $R^{5}_{q}$, it is just a matter of changing $\sum_{q}$ into
$\sup_{q}$. We proceed similarly for $R^{4}_{q}$.
\null{\hfill $\Box$}


\begin{thebibliography}{}
\bibitem{Ab}
H. ABIDI. \'Equation de Navier-Stokes avec densit\'e et viscosit\'e variables dans l'espace critique.  \textit{Th\`{e}se de l'Universit\'e Paris VI}.
\bibitem{AP}
H. ABIDI and M. PAICU. \'Equation de Navier-Stokes avec densit\'e et viscosit\'e variables dans l'espace critique. \textit{Annales de l'institut Fourier}, 57 no. 3 (2007), p. 883-917.
\bibitem{AKM}
S. ANTONTSEV, A. KAZHIKOV and V. MONAKOV. Boundary value problems in mechanics of nonhomogeneous fluids. Translated from the
Russian.  \textit{Studies in mathematics and its applications, 22. North-Holland publishing co. Amsterdam}, 1990.
\bibitem{BC}
H. BAHOURI and J.-Y. CHEMIN. \'Equations d'ondes quasilin\'eaires et
estimation de Strichartz,  \textit{Amer. J. Mathematics 121}, (1999), 1337-1377.
\bibitem{BCD}
H. BAHOURI, J.-Y. CHEMIN and R. DANCHIN. Fourier analysis and nonlinear partial differential equations,
\textit{to appear in Springer}.
\bibitem{Bo}
J.-M. BONY. Calcul symbolique et propagation des singularit\'es pour
les \'equations aux d\'eriv\'ees partielles non lin\'eaires, \textit{Annales
Scientifiques de l'\'ecole Normale Sup\'erieure 14}, (1981),
209-246.
\bibitem{CMP}
M. CANNONE, Y. MEYER and F. PLANCHON. Solutions auto-similaires des \'equations de Navier-Stokes. \textit{S\'eminaire sur les \'equations
aux d\'eriv\'ees partielles, 1993-1994}, exp. No12 pp. \'Ecole polytech, palaiseau, 1994.
\bibitem{Ch1}
J-Y CHEMIN. Fluides parfaits incompressibles. \textit{Ast\'erisque, 230}, 1995.
\bibitem{Ch2}
J.-Y. CHEMIN. Th\'eor\`emes d'unicit\'e pour le syst\`eme de
Navier-Stokes tridimensionnel, \textit{J. d'Analyse Math.}, 77, (1999), 25-50.
\bibitem{Ch3}
J.-Y. CHEMIN and N. LERNER. Flot de champs de vecteurs non
lipschitziens et \'equations de Navier-Stokes, \textit{J. Differential
Equations.}, 121, (1995), 247-286.
\bibitem{Da1}
R. DANCHIN. Fourier analysis method for PDE's, \textit{Preprint Novembre 2005}.
\bibitem{Da2}
R. DANCHIN. Local and global well-posedness results for flows of inhomogeneous viscous fluids. \textit{Advances in differential
equations}, 9 (2004), 353-386.
\bibitem{Da3}
R. DANCHIN. Density-dependent incompressible viscous fluids in critical spaces. \textit{Proc. Roy. Soc. Edinburgh Sect. A}, 133(6),
1311-1334, 2003.
\bibitem{Da4}
R. DANCHIN. The inviscid limit for density-dependent incompressible fluids. \textit{Ann.
Fac. Sci. Toulouse Math.}, S\'er. 6, 15 no. 4 (2006), p. 637-688.
\bibitem{Da5}
R. DANCHIN. Well-posedness in critical spaces for barotropic viscous fluids with truly not constant density. \textit{Communications in Partial Differential Equations},32:9,1373 — 1397, 2007.
\bibitem{Da6}
R. DANCHIN. Local Theory in critical Spaces for Compressible Viscous
and Heat-Conductive Gases,[ \textit{Communication in
Partial Differential Equations.} 26 (2001), 7-8, 1183-1233]. \textit{Commun. Partial differential equations.} 27, (2002), 11-12, 2531-2532.
\bibitem{Da7}
R. DANCHIN. Global Existence in Critical Spaces for Flows of
Compressible Viscous and Heat-Conductive Gases, \textit{Arch. Rational
Mech. Anal.}, 160, (2001), 1-39.
\bibitem{Da8}
R. DANCHIN. On the uniqueness in critical spaces for compressible Navier-Stokes equations. \textit{NoDEA Nonlinear Differentiel Equations
Appl}, 12(1), 111-128, 2005.
\bibitem{De1}
B. DESJARDINS. Linear transport equations with initial values in Sobolev spaces and application to the Navier-Stokes
equations, \textit{Differential and Integral Equations}, 10 (1997), 587-598.
\bibitem{De2}
B. DESJARDINS. Global existence results for the incompressible density-dependent Navier-Stokes equations in the whole
space. \textit{Differential and integral equations}, 10 (1997), 587-598.
\bibitem{De3}
B. DESJARDINS. Regularity results for two-dimensional flows of multiphase viscous fluids. \textit{Arch. Rational mech. Anal.}, 137, (1997),
135-158.
\bibitem{FCG}
E. FERNANDEZ-CARA and F. GUILLEN. The existence of nonhomogeneous viscous and incompressible flow in unbounded domains. \textit{Comm.
in Partial Differential Equation}, 17, (1992), 1253-1265.
\bibitem{FK}
H. FUJITA and T. KATO. On the Navier-Stokes initial value problem I, \textit{Archive for Rational Mechanics and Analysis},
16, (1964), 269-315.
\bibitem{IT}
S. ITOH and A. TANI. Solvability of nonstationnary problems for nonhomogeneous
incompressible fluids and the convergence with vanishing viscosity, \textit{Tokyo Journal
of Mathematics}, 22, (1999), 17-42.
\bibitem{Ka}
V. KAZHIKOV. Resolution of boundary value problems for nonhomogeneous viscous fluids. \textit{Dokl. Akad. Nauh}, 216, (1974),
1008-1010.
\bibitem{KT}
H. KOCH and D. TATARU. Well-posedness for the Navier-Stokes equations, \textit{Advances in Math}. 157, 920010, 22-35.
\bibitem{LS}
O. LADYZHENSKAYA and V. SOLONNIKOV. The unique solvability of an initial-boundary value problem for viscous incompressible
inhomogeneous fluids. \textit{J. Soviet Math.} 9, (1978), 697-749.
\bibitem{Li}
P-. L. LIONS, Mathematical topics in fluid dynamics. Vol 1 incompressible models, \textit{Oxford university press} 1996.
\bibitem{Me}
Y. MEYER. Wavelets, paraproducts, and Navier-Stokes equation.
In Current developments in mathematics, 1996 (Cambridge, MA),
page 105-212. \textit{Int. Press, Boston}, MA, 1997.
\bibitem{RS}
T. RUNST and W. SICKEL, Sobolev spaces of fractional order, Nemytskij
operators, and nonlinear partial differential equations. \textit{de Gruyter Series in
Nonlinear Analysis and Applications, 3. Walter de Gruyter and Co.}, Berlin (1996)\\
\end{thebibliography}
\end{document}